\documentclass[a4paper,12pt,reqno]{amsart}

\usepackage{amsfonts, color}
\usepackage{amsmath}
\usepackage{amssymb}
\usepackage[a4paper]{geometry}
\usepackage{mathrsfs}
\usepackage[colorlinks]{hyperref}
\usepackage{graphicx}

\renewcommand\eqref[1]{(\ref{#1})}

\makeatletter
\newcommand*{\mint}[1]{%
  \mint@l{#1}{}%
}
\newcommand*{\mint@l}[2]{%
  \@ifnextchar\limits{%
    \mint@l{#1}%
  }{%
    \@ifnextchar\nolimits{%
      \mint@l{#1}%
    }{%
      \@ifnextchar\displaylimits{%
        \mint@l{#1}%
      }{%
        \mint@s{#2}{#1}%
      }%
    }%
  }%
}
\newcommand*{\mint@s}[2]{%
  \@ifnextchar_{%
    \mint@sub{#1}{#2}%
  }{%
    \@ifnextchar^{%
      \mint@sup{#1}{#2}%
    }{%
      \mint@{#1}{#2}{}{}%
    }%
  }%
}
\def\mint@sub#1#2_#3{%
  \@ifnextchar^{%
    \mint@sub@sup{#1}{#2}{#3}%
  }{%
    \mint@{#1}{#2}{#3}{}%
  }%
}
\def\mint@sup#1#2^#3{%
  \@ifnextchar_{%
    \mint@sup@sub{#1}{#2}{#3}%
  }{%
    \mint@{#1}{#2}{}{#3}%
  }%
}
\def\mint@sub@sup#1#2#3^#4{%
  \mint@{#1}{#2}{#3}{#4}%
}
\def\mint@sup@sub#1#2#3_#4{%
  \mint@{#1}{#2}{#4}{#3}%
}
\newcommand*{\mint@}[4]{%
  \mathop{}%
  \mkern-\thinmuskip
  \mathchoice{%
    \mint@@{#1}{#2}{#3}{#4}%
        \displaystyle\textstyle\scriptstyle
  }{%
    \mint@@{#1}{#2}{#3}{#4}%
        \textstyle\scriptstyle\scriptstyle
  }{%
    \mint@@{#1}{#2}{#3}{#4}%
        \scriptstyle\scriptscriptstyle\scriptscriptstyle
  }{%
    \mint@@{#1}{#2}{#3}{#4}%
        \scriptscriptstyle\scriptscriptstyle\scriptscriptstyle
  }%
  \mkern-\thinmuskip
  \int#1%
  \ifx\\#3\\\else_{#3}\fi
  \ifx\\#4\\\else^{#4}\fi
}
\newcommand*{\mint@@}[7]{%
  \begingroup
    \sbox0{$#5\int\m@th$}%
    \sbox2{$#5\int_{}\m@th$}%
    \dimen2=\wd0 %
    \let\mint@limits=#1\relax
    \ifx\mint@limits\relax
      \sbox4{$#5\int_{\kern1sp}^{\kern1sp}\m@th$}%
      \ifdim\wd4>\wd2 %
        \let\mint@limits=\nolimits
      \else
        \let\mint@limits=\limits
      \fi
    \fi
    \ifx\mint@limits\displaylimits
      \ifx#5\displaystyle
        \let\mint@limits=\limits
      \fi
    \fi
    \ifx\mint@limits\limits
      \sbox0{$#7#3\m@th$}%
      \sbox2{$#7#4\m@th$}%
      \ifdim\wd0>\dimen2 %
        \dimen2=\wd0 %
      \fi
      \ifdim\wd2>\dimen2 %
        \dimen2=\wd2 %
      \fi
    \fi
    \rlap{%
      $#5%
        \vcenter{%
          \hbox to\dimen2{%
            \hss
            $#6{#2}\m@th$%
            \hss
          }%
        }%
      $%
    }%
  \endgroup
}

\numberwithin{equation}{section}
\theoremstyle{plain}
\newtheorem{thm}{Theorem}[section]
\newtheorem{lem}[thm]{Lemma}
\newtheorem{prop}[thm]{Proposition}

\theoremstyle{definition}

\theoremstyle{remark}
\newtheorem{rem}{Remark}[section]
\newtheorem{defn}{Definition}

\numberwithin{equation}{section}

\begin{document}

\title[Parabolic Equations with Singular Coefficients and Boundary Data]{Parabolic Equations with Singular Coefficients and Boundary Data: Analysis and Numerical Simulations}

\author[A. Altybay]{Arshyn Altybay}
\address{
  Arshyn Altybay:
  \endgraf
  Department of Mathematics: Analysis, Logic and Discrete Mathematics
  \endgraf
  Ghent University, Belgium
  \endgraf
  and
  \endgraf
  Institute of Mathematics and Mathematical Modelling
  \endgraf
  125 Pushkin str., 050010 Almaty, Kazakhstan
  \endgraf
   and
  \endgraf 
    al--Farabi Kazakh National University
  \endgraf
  71 al--Farabi ave., 050040 Almaty, Kazakhstan
  \endgraf
  {\it E-mail address} {\rm arshyn.altybay@gmail.com}
}

\author[A. Yeskermessuly ]{Alibek Yeskermessuly}
\address{
  Alibek Yeskermessuly:
  \endgraf   
  Faculty of Natural Sciences and Informatization
  \endgraf
  Altynsarin University, Auelbekov, 17, 110300 Arkalyk, Kazakhstan
  \endgraf  
  {\it E-mail address:} {\rm alibek.yeskermessuly@gmail.com}
  }

\thanks{This research was funded by the Science Committee of the Ministry of Science and Higher Education of the Republic of Kazakhstan (Grant No. AP23486342), by the FWO Research Grant G083525N: Evolutionary partial differential equations with strong singularities, and by the Methusalem programme of the Ghent University Special Research Fund (BOF) (Grant number 01M01021).}

\keywords{parabolic equation, divergence form, drift term, energy methods, Galerkin approximation, singular coefficients, regularisation, very weak solution.}
\subjclass[2020]{35K20, 35D30, 35K67.}

\begin{abstract}
We investigate linear parabolic equations in divergence form with singular coefficients and nonsmooth initial--boundary data. When the diffusion, drift, or potential terms, as well as the source term and boundary conditions, are distributions rather than functions, classical and weak solution concepts become inadequate, since products involving distributions are not well defined in general. To address this difficulty, we introduce a framework of very weak solutions based on regularisation procedures and the theory of moderate nets. Under the stated moderateness and uniform-ellipticity assumptions, we establish existence of very weak solutions and prove uniqueness via negligibility arguments. Moreover, in the regular-data regime, we show consistency with classical weak solutions. Finally, we present numerical experiments illustrating the behaviour of regularised solutions for highly singular inputs, including delta-type potentials and distributional boundary traces.
\end{abstract}

\maketitle

\section{Introduction}

In this paper, we study a linear parabolic equation in divergence form with a drift term and potential on a bounded one-dimensional domain. We focus on the setting where the coefficients and the initial-boundary data are allowed to be singular, in the sense of distributions. The diffusion coefficient may also be singular, but it is required to be uniformly positive in the distributional sense described below.

More precisely, we consider the initial--boundary value problem
\begin{align}\label{1.1}
    \partial_t u(t,x) - \partial_x \left(a(t,x)\partial_x u(t,x)\right) 
    + b(t,x) \partial_x u(t,x) 
    + q(x) u(t,x) 
    = f(t,x),\nonumber\\
    \quad (t,x) \in [0,T] \times (0,1),
\end{align}
with initial condition
\begin{equation}\label{1.2}
    u(0,x) = u_0(x), \quad x \in (0,1),
\end{equation}
and inhomogeneous Dirichlet boundary conditions
\begin{equation}\label{1.3}
    u(t,0) = g_0(t), \quad u(t,1) = g_1(t), \quad t \in [0,T].
\end{equation}
Here \(a(t,x)\) denotes the diffusion coefficient, \(b(t,x)\) the drift field,
\(q(x)\) a potential term, and \(f(t,x)\) an external source.
The functions \(u_0\), \(g_0\), and \(g_1\) prescribe the initial and boundary data.

The diffusion coefficient cannot be an arbitrary distribution. Throughout the singular theory, we impose the distributional ellipticity condition
\[
\langle a-\alpha,\varphi\rangle\ge0
\quad\text{for every }0\le\varphi\in C_c^\infty((0,T)\times(0,1)),
\]
for some \(\alpha>0\). Equivalently, \(a=\alpha+\mu\), where \(\mu\) is a nonnegative Radon measure. Thus, singularities in the diffusivity are permitted only through a nonnegative measure component superimposed on a strictly positive background. This is precisely the class for which nonnegative mollification produces regularised coefficients satisfying \(a_\varepsilon\ge\alpha\) uniformly in \(\varepsilon\).

The variational well-posedness theory for linear parabolic equations
with regular coefficients and boundary data is classical; see, for
example, \cite[Section~7.1]{Evans} and \cite{LionsMagenes}. The
difficulty considered here arises when one or more coefficients or data
are distributions. In that case, products such as
\(b\,\partial_xu\) and \(qu\), and more generally the interaction
between distributional coefficients and the unknown solution, are not
defined within the classical distributional framework, in accordance
with Schwartz's impossibility theorem \cite{Sch54}.

A systematic framework for nonlinear operations involving
distributions is provided by Colombeau-type generalised functions; see
the monograph \cite{Ober}. In the context of partial differential
equations, Garetto and Ruzhansky \cite{GR15} introduced the notion of a
very weak solution for hyperbolic equations with non-regular
coefficients. Subsequent developments treated wave and hyperbolic
equations with irregular propagation speeds, electromagnetic fields,
or space-dependent singular coefficients
\cite{RT17a, RT17b, Gar21}. For heat-type equations, strongly singular
potentials, hypoelliptic operators, and irregular thermal conductivity
with distributional data were considered in
\cite{ARST21c, CDRT23, CRT22a}, respectively.

For evolution equations associated with Sturm--Liouville operators,
very weak solutions have been constructed by means of spectral
decompositions and separation of variables. In particular,
\cite{RSY22, RY22, ARST25} study wave equations with singular coefficients,
while \cite{RY24a} treats the corresponding heat equation. These
results concern operator settings in which the spectral structure
provides a diagonal representation of the regularised problems.

Parabolic equations with low-regularity coefficients or data have also
been studied using solution concepts different from the very weak
solution framework. Issoglio \cite{Issoglio19} considers a nonlinear
parabolic equation with a distributional coefficient belonging to a
negative-order Besov space and establishes well-posedness in a mild
solution framework. Pierre and Vovelle \cite{PierreVovelle12} use a
kinetic formulation to obtain renormalized solutions for a nonlinear
parabolic Cauchy--Dirichlet problem with \(L^1\)-data. Michel and
Vovelle \cite{MichelVovelle03} develop an entropy formulation for
degenerate parabolic--hyperbolic problems with general nonhomogeneous
Dirichlet boundary conditions and prove convergence of finite-volume
approximations. Di~Nardo, Feo, and Guibé
\cite{DiNardoFeoGuibe13} establish uniqueness of renormalized
solutions for nonlinear parabolic problems with low-regularity
right-hand sides and lower-order terms. These approaches use mild,
kinetic, entropy, or renormalized solution concepts and therefore
differ from the regularisation-based very weak solution framework
adopted in the present paper.

The present contribution differs from the preceding works in several
respects. First, the diffusion coefficient is allowed to be
measure-valued,
\[
a=\alpha+\mu,
\qquad
\alpha>0,
\qquad
\mu\in\mathcal M_+(Q_T),
\]
and is regularised by a positivity-preserving procedure that retains
the uniform bound \(a_\varepsilon\geq\alpha\). Second, the drift,
potential, source, initial datum, and nonhomogeneous Dirichlet boundary
data may be distributional simultaneously. Third, instead of relying
on a spectral diagonalisation, we use coefficient-dependent
variational energy estimates for the regularised divergence-form
problems. These estimates yield moderateness of the solution net,
uniqueness modulo negligible nets, and consistency with the classical
weak solution whenever the coefficients and data possess the required
regularity.

Singular coefficients and data model abrupt interfaces, localised
inclusions, point interactions, and concentrated sources in, for example,
heterogeneous diffusion, heat conduction, and drift--diffusion problems. Such
idealised data may fall outside the regularity required by classical PDE
theory, which motivates a regularisation-based solution concept and a careful
interpretation of the resulting solution net.

The framework provides a basis for numerical approximation because, for each fixed \(\varepsilon>0\), it yields a regular problem to which standard discretisation methods apply. Moderateness controls the growth of this family, but does not by itself imply convergence as \(\varepsilon\to0\); convergence is obtained in the regular-data regime through the consistency result below.

Accordingly, we prove existence and uniqueness of very weak solutions to
\eqref{1.1}--\eqref{1.3}, establish consistency in the regular-data regime, and
illustrate the behaviour of regularised solutions for delta-type inputs.

The present work continues our recent analysis of the Klein--Gordon equation with singular coefficients
under inhomogeneous Dirichlet boundary conditions \cite{RY25}.
Many of the analytical techniques developed there, including the regularisation strategy
and the derivation of uniform \emph{a priori} estimates, provide the foundation for the current study.
In contrast to the hyperbolic setting, the parabolic case introduces additional difficulties,
since singularities interact with the time evolution and require new arguments.

To the best of our knowledge, well-posedness for parabolic divergence-form equations
with distributional coefficients and nonhomogeneous Dirichlet boundary data
has not previously been established within the very weak solution framework.
Thus, the primary novelty of this work lies in the analysis of a parabolic equation
with singular diffusion, drift, and potential terms, together with a singular source term,
under inhomogeneous Dirichlet boundary conditions.

Although the analytical results of this paper are formulated on the interval
\((0,1)\), the regularisation strategy is not intrinsically one-dimensional. On a
bounded Lipschitz domain \(\Omega\subset\mathbb{R}^{d}\), the natural analogue is
\[
\partial_t u_\varepsilon
-\operatorname{div}(A_\varepsilon\nabla u_\varepsilon)
+b_\varepsilon\cdot\nabla u_\varepsilon
+q_\varepsilon u_\varepsilon
=f_\varepsilon,
\]
where the regularised diffusion matrices satisfy
\[
A_\varepsilon(t,x)\xi\cdot\xi
\geq \alpha|\xi|^2
\]
with a constant \(\alpha>0\) independent of \(\varepsilon\). The basic
variational and energy arguments can then be formulated in
\(H_0^1(\Omega)\), while inhomogeneous boundary data can be treated by suitable
trace and lifting operators. A complete extension would nevertheless require
additional assumptions on the domain and on matrix-valued distributional
coefficients.

The same philosophy may also be applied to certain semilinear or quasilinear
equations, provided that the nonlinearities satisfy structural assumptions, such
as Lipschitz continuity, monotonicity, or controlled growth, which yield
\(\varepsilon\)-uniform a priori estimates. The treatment of general nonlinear
products of distributions lies beyond the scope of the present paper.

The paper is organised as follows.
In Section~2 we recall the energy estimates for the homogeneous Dirichlet problem.
Section~3 contains the existence, uniqueness, and consistency results for very weak solutions
under inhomogeneous Dirichlet boundary conditions.
Section~4 is devoted to numerical experiments.

\section{The regular problem and coefficient-dependent estimates}
\label{sec:regular-problem}

Let
\[
Q_T:=(0,T)\times(0,1).
\]

We consider a linear parabolic equation in divergence form with drift and potential terms, subject to an initial condition and homogeneous Dirichlet boundary conditions:
\begin{equation}\label{hom.problem}
\left\{\begin{array}{l}
\partial_t u - \partial_x(a(t,x)\partial_x u) + b(t,x) \partial_x u + q(x) u= f(t,x), \quad (t,x)\in Q_T,\\
u(0,x)=u_0(x),\quad x\in (0,1),\\
u(t,0)=u(t,1)=0, \quad t\in (0,T).
\end{array}\right.
\end{equation}
Here \(a,b\in L^\infty(Q_T)\) and
\(q\in L^\infty(0,1)\) are real-valued coefficients,
\[
a(t,x)\geq\alpha>0
\quad\text{for almost every }(t,x)\in Q_T,
\]
while
\[
f\in L^2(0,T;L^2(0,1)),
\qquad
u_0\in L^2(0,1).
\]
We allow \(q\) to change sign and write
\[
q=q_+-q_-,
\qquad
q_+:=\max\{q,0\},
\qquad
q_-:=\max\{-q,0\}.
\]

The well-posedness theory for linear parabolic equations with bounded coefficients is classical; see, for example, \cite[Section~7.1]{Evans} and \cite{LionsMagenes}. The Galerkin method used in the derivation of the estimates below is also standard; see, for instance, \cite{Thomee06}. We therefore omit the standard Galerkin existence proof, including the compactness, passage-to-the-limit, and uniqueness arguments, and retain only the coefficient-dependent estimates needed in the subsequent analysis of very weak solutions. In particular, we keep track explicitly of their dependence on the ellipticity constant and on the coefficient norms, since these quantities depend on the regularisation parameter in Section~\ref{sec:very-weak}.

For notational convenience, the variational formulation is written in the complex Hilbert spaces \(L^2(0,1;\mathbb C)\) and \(H_0^1(0,1;\mathbb C)\), endowed with the standard sesquilinear
inner product
\[
(u,v)_{L^2(0,1)}
:=\int_0^1u(x)\overline{v(x)}\,dx.
\]

We use the following notion of weak solution.
\begin{defn}[Weak solution]\label{def:weak_solution}
A function
\[
u\in L^2(0,T;H_0^1(0,1))
\cap C([0,T];L^2(0,1)),
\qquad
\partial_tu\in L^2(0,T;H^{-1}(0,1)),
\]
is called a \emph{weak solution} of \eqref{hom.problem} if
\(H^{-1}(0,1)\) is the dual of \(H_0^1(0,1)\),
\(u(0)=u_0\) in \(L^2(0,1)\) and, for every
\(\varphi\in H_0^1(0,1)\) and almost every \(t\in(0,T)\),
\begin{align}\label{eq:weak_formulation}
\langle \partial_t u(t),\varphi\rangle_{H^{-1},H^1_0}
&+ \int_0^1 a(t,x)\, \partial_x u(t,x)\,\overline{\partial_x \varphi(x)} \, dx
+ \int_0^1 b(t,x)\, \partial_x u(t,x)\, \overline{\varphi(x)} \, dx \nonumber\\
&+ \int_0^1 q(x)\, u(t,x)\, \overline{\varphi(x)} \, dx
= \int_0^1 f(t,x)\, \overline{\varphi(x)}\, dx.
\end{align}
\end{defn}

\subsection{Galerkin approximation}

We briefly recall the Galerkin approximation used to derive the coefficient-dependent estimates below. Let
\[
w_k(x):=\sqrt{2}\sin(k\pi x),
\qquad k\in\mathbb N.
\]
Then \(\{w_k\}_{k=1}^{\infty}\) is an orthonormal basis of \(L^2(0,1)\) and an orthogonal basis of \(H_0^1(0,1)\).
With
\[
(u,v)_{H_0^1(0,1)}:=\int_0^1u_x(x)\overline{v_x(x)}\,dx,
\]
the basis relations are
\[
(w_j,w_k)_{L^2(0,1)}=\delta_{jk},
\qquad
(w_j,w_k)_{H_0^1(0,1)}=(k\pi)^2\delta_{jk}.
\]

For each \(m\in\mathbb N\), set
\[
V_m:=\operatorname{span}\{w_1,\ldots,w_m\},
\]
and let
\[
P_m:L^2(0,1)\to V_m
\]
denote the \(L^2(0,1)\)-orthogonal projection.

We seek a Galerkin approximation \(u_m\in H^1(0,T;V_m)\) of the form
\begin{equation}\label{u_m}
u_m(t,x)=\sum_{k=1}^{m}d_m^k(t)w_k(x),
\end{equation}
where
\[
u_m(0)=P_mu_0,
\qquad\text{equivalently}\qquad
d_m^k(0)=(u_0,w_k)_{L^2(0,1)},
\quad k=1,\ldots,m,
\]
and
\begin{equation}\label{u_m:weak_form}
(\partial_tu_m,w_k)
+(a\partial_xu_m,\partial_xw_k)
+(b\partial_xu_m,w_k)
+(qu_m,w_k)
=(f,w_k),
\qquad k=1,\ldots,m,
\end{equation}
for almost every \(t\in(0,T)\).

This yields a finite-dimensional linear system with measurable bounded coefficients and an \(L^2\)-right-hand side. Standard linear ODE theory therefore gives a unique Galerkin solution
\[
u_m\in H^1(0,T;V_m)\subset C([0,T];V_m).
\]

\subsection{Energy estimates}
\begin{thm}\label{thm:energy-estimate}
Let $u_m$ be the Galerkin approximation defined in \eqref{u_m}, and assume that the coefficients satisfy
\begin{equation}\label{coefficients}
a,b,\partial_xb\in L^\infty(Q_T),\quad
a(t,x)\geq\alpha>0
\text{  for a.e. }(t,x)\in Q_T,
\quad
q\in L^\infty(0,1).
\end{equation}
Then, for all \( m \in \mathbb{N} \), the following energy estimate holds:
\begin{align}\label{eq:energy-estimate}
&\| u_m \|_{L^\infty(0,T; L^2(0,1))} 
+ \sqrt{\alpha}\| \partial_x u_m \|_{L^2(0,T; L^2(0,1))} 
+ \| \sqrt{q_+} u_m \|_{L^2(0,T; L^2(0,1))} \nonumber\\
&\qquad+ \| \partial_t u_m \|_{L^2(0,T; H^{-1}(0,1))} 
\leq C \cdot \exp{\left\{T\big(\|\partial_x b\|_{L^\infty(Q_T)}+2\|q_-\|_{L^\infty(0,1)}\big)\right\}}\nonumber\\
&\times \left(1 + \frac{1}{\sqrt{\alpha}}\left(\| a \|_{L^\infty(Q_T)} + \| b \|_{L^\infty(Q_T)}\right) + \|q \|_{L^\infty(0,1)} \right)\nonumber \\
&\times\left( \| u_0 \|_{L^2(0,1)} + \| f \|_{L^2(0,T; L^2(0,1))} \right),
\end{align}
where $C>0$ depends only on $T$ and the spatial domain $(0,1)$, and is independent of $m$.
All dependence on the coefficients $a,b,q$ (and on the data $u_0,f$) is explicitly displayed on the right-hand side of \eqref{eq:energy-estimate}.
\end{thm}

\begin{proof}
By linearity, the Galerkin identity \eqref{u_m:weak_form}
holds for every test function \(v_m\in V_m\). Taking
\(v_m=u_m(t)\) and then taking real parts, we obtain
\[
\frac12\frac{d}{dt}\|u_m\|_{L^2(0,1)}^2
+(a\partial_xu_m,\partial_xu_m)
+\operatorname{Re}(b\partial_xu_m,u_m)
+\int_0^1q(x)|u_m(t,x)|^2\,dx
=\operatorname{Re}(f,u_m).
\]
Since the coefficients are real-valued and \(u_m\) satisfies the
homogeneous Dirichlet condition,
\[
\operatorname{Re}(b\partial_xu_m,u_m)
=-\frac12\int_0^1
\partial_xb(t,x)|u_m(t,x)|^2\,dx.
\]
Using \(a\geq\alpha\), \(q=q_+-q_-\), and Young's inequality, we
therefore obtain
\begin{align}\label{eq:differential-energy}
\frac{d}{dt}\|u_m\|_{L^2(0,1)}^2
&+2\alpha\|\partial_xu_m\|_{L^2(0,1)}^2
+2\|\sqrt{q_+}u_m\|_{L^2(0,1)}^2
\nonumber\\
&\leq
\|f\|_{L^2(0,1)}^2
+C_1\|u_m\|_{L^2(0,1)}^2,
\end{align}
where
\[
C_1
:=
1+\|\partial_xb\|_{L^\infty(Q_T)}
+2\|q_-\|_{L^\infty(0,1)}.
\]

Set
\[
B:=\|\partial_xb\|_{L^\infty(Q_T)}
+2\|q_-\|_{L^\infty(0,1)},
\]
so that \(C_1=1+B\). Multiplying
\eqref{eq:differential-energy} by \(e^{-C_1t}\) and integrating,
we obtain
\begin{align*}
&\|u_m\|_{L^\infty(0,T;L^2(0,1))}
+\sqrt{\alpha}\,
 \|\partial_xu_m\|_{L^2(0,T;L^2(0,1))}\\
&\qquad
+\|\sqrt{q_+}\,u_m\|_{L^2(0,T;L^2(0,1))}
\\
&\leq
C(T)e^{BT/2}
\left(
\|u_0\|_{L^2(0,1)}
+\|f\|_{L^2(0,T;L^2(0,1))}
\right),
\end{align*}
where we used \(u_m(0)=P_mu_0\),
\(\|P_mu_0\|_{L^2}\leq\|u_0\|_{L^2}\), and absorbed
\(e^{T/2}\) into \(C(T)\). Since \(B\geq0\), we have
\[
e^{BT/2}\leq e^{BT}.
\]
Therefore, the preceding estimate proves the bounds for the first three terms on the left-hand side of
\eqref{eq:energy-estimate}, with the exponential factor displayed there.

Since \(\partial_tu_m(t)\in V_m\), for every
\(\phi\in H_0^1(0,1)\) we have
\[
\langle\partial_tu_m(t),\phi\rangle
=
(\partial_tu_m(t),P_m\phi)_{L^2(0,1)}.
\]
Using \(P_m\phi\in V_m\) as a test function in the Galerkin
identity gives
\begin{align*}
|\langle\partial_tu_m(t),\phi\rangle|
&\leq
\|f(t)\|_{L^2(0,1)}
\|P_m\phi\|_{L^2(0,1)}\\
&\quad+
\|a(t,\cdot)\|_{L^\infty(0,1)}
\|\partial_xu_m(t)\|_{L^2(0,1)}
\|\partial_xP_m\phi\|_{L^2(0,1)}\\
&\quad+
\|b(t,\cdot)\|_{L^\infty(0,1)}
\|\partial_xu_m(t)\|_{L^2(0,1)}
\|P_m\phi\|_{L^2(0,1)}\\
&\quad+
\|q\|_{L^\infty(0,1)}
\|u_m(t)\|_{L^2(0,1)}
\|P_m\phi\|_{L^2(0,1)}.
\end{align*}
For the sine basis, \(P_m\) is stable both in \(L^2(0,1)\)
and in \(H_0^1(0,1)\):
\[
\|P_m\phi\|_{L^2(0,1)}
\leq\|\phi\|_{L^2(0,1)},
\qquad
\|\partial_xP_m\phi\|_{L^2(0,1)}
\leq\|\partial_x\phi\|_{L^2(0,1)}.
\]
Consequently,
\begin{align*}
\|\partial_tu_m\|_{L^2(0,T;H^{-1}(0,1))}
&\leq C\Big[
\|f\|_{L^2(0,T;L^2(0,1))}\\
&\quad+
\big(
\|a\|_{L^\infty(Q_T)}
+\|b\|_{L^\infty(Q_T)}
\big)
\|\partial_xu_m\|_{L^2(0,T;L^2(0,1))}\\
&\quad+
\|q\|_{L^\infty(0,1)}
\|u_m\|_{L^2(0,T;L^2(0,1))}
\Big].
\end{align*}
Using the preceding bounds and
\[
\|u_m\|_{L^2(0,T;L^2(0,1))}
\leq
\sqrt{T}\,
\|u_m\|_{L^\infty(0,T;L^2(0,1))},
\]
we obtain
\begin{align*}
\|\partial_tu_m\|_{L^2(0,T;H^{-1}(0,1))}
&\leq
C
\exp\!\left\{
T\left(
\|\partial_xb\|_{L^\infty(Q_T)}
+2\|q_-\|_{L^\infty(0,1)}
\right)
\right\}\\
&\quad\times
\left(
1+
\frac{
\|a\|_{L^\infty(Q_T)}
+\|b\|_{L^\infty(Q_T)}
}{\sqrt{\alpha}}
+\|q\|_{L^\infty(0,1)}
\right)\\
&\quad\times
\left(
\|u_0\|_{L^2(0,1)}
+\|f\|_{L^2(0,T;L^2(0,1))}
\right).
\end{align*}
Combining this estimate with the first three bounds proves
\eqref{eq:energy-estimate}.
\end{proof}

\begin{thm}\label{thm:1D-wellposedness}
Assume that
\[
a,b,\partial_xb\in L^\infty(Q_T),\qquad
a(t,x)\geq\alpha>0
\quad\text{for a.e. }(t,x)\in Q_T,
\qquad
q\in L^\infty(0,1),
\]
and let
\[
u_0\in L^2(0,1),
\qquad
f\in L^2(0,T;L^2(0,1)).
\]
Then problem \eqref{hom.problem} admits a unique weak solution
\[
u\in L^2(0,T;H_0^1(0,1))
\cap C([0,T];L^2(0,1)),
\qquad
\partial_tu\in L^2(0,T;H^{-1}(0,1)).
\]
Moreover,
\[
\sqrt{q_+}\,u\in L^2(0,T;L^2(0,1)),
\]
and estimate \eqref{eq:energy-estimate}, with \(u_m\) replaced by
\(u\), holds for the weak solution.
\end{thm}
\begin{proof}
For almost every \(t\in(0,T)\), define
\[
\mathcal B_t(u,v)
:=
\int_0^1 a(t,x)u_x\overline{v_x}\,dx
+\int_0^1 b(t,x)u_x\overline{v}\,dx
+\int_0^1 q(x)u\overline{v}\,dx.
\]
For every \(u,v\in H_0^1(0,1)\), the map
\(t\mapsto\mathcal B_t(u,v)\) is measurable, and
\(\mathcal B_t\) is uniformly bounded on
\(H_0^1(0,1)\times H_0^1(0,1)\). Moreover,
\[
\operatorname{Re}\mathcal B_t(v,v)
\geq
\alpha\|v_x\|_{L^2(0,1)}^2
-
\left(
\frac12\|\partial_xb\|_{L^\infty(Q_T)}
+\|q_-\|_{L^\infty(0,1)}
\right)
\|v\|_{L^2(0,1)}^2.
\]
Thus, \(\mathcal B_t\) satisfies the standard G{\aa}rding inequality.
Existence, uniqueness, and time continuity follow from the standard
variational theory for linear parabolic equations; see, for example,
\cite[Section~7.1]{Evans} and
\cite{LionsMagenes}.

Applying Theorem~\ref{thm:energy-estimate} to the Galerkin approximations and passing to the limit in the standard variational construction, using weak lower semicontinuity, yields the stated coefficient-dependent estimate for \(u\).
\end{proof}

The following higher-regularity estimate is retained because its explicit dependence on the coefficient norms is used in the moderateness and consistency arguments of Section~\ref{sec:very-weak}.

\begin{lem}[Higher time--space regularity]
\label{lem:time-space-second-derivative}
Let \(u\) be the unique weak solution of \eqref{hom.problem}
given by Theorem~\ref{thm:1D-wellposedness}. Assume that
\[
u_0\in H_0^1(0,1),
\qquad
f\in L^2(0,T;L^2(0,1)),
\]
and that the real-valued coefficients satisfy
\[
a\in W^{1,\infty}(Q_T),
\qquad
b\in L^\infty(0,T;W^{1,\infty}(0,1)),\qquad q\in L^\infty(0,1),
\]
and
\[
a(t,x)\geq\alpha>0 \quad \text{for almost every } (t,x)\in Q_T.
\]
Then
\[
\partial_tu\in L^2(0,T;L^2(0,1)),
\qquad
\partial_xu\in L^\infty(0,T;L^2(0,1)),
\]
and
\[
\partial_x(a\partial_xu),\;
\partial_x^2u
\in L^2(0,T;L^2(0,1)).
\]

Moreover,
\begin{align}
&\|\partial_tu\|_{L^2(0,T;L^2(0,1))}
+\sqrt{\alpha}\,
 \|\partial_xu\|_{L^\infty(0,T;L^2(0,1))}
+\|\partial_x(a\partial_xu)\|_{L^2(0,T;L^2(0,1))}
\nonumber\\
&\qquad
+\alpha\|\partial_x^2u\|_{L^2(0,T;L^2(0,1))}
\leq
C\exp\!\left\{
T\left(
\|\partial_xb\|_{L^\infty(Q_T)}
+2\|q_-\|_{L^\infty(0,1)}
\right)
\right\}
\nonumber\\
&\qquad\quad\times
\bigg[
1+\|a\|_{L^\infty(Q_T)}^{1/2}
+\|q\|_{L^\infty(0,1)}
\nonumber\\
&\hspace{3cm}
+\frac{1}{\sqrt{\alpha}}
\left(
\|\partial_ta\|_{L^\infty(Q_T)}^{1/2}
+\|\partial_xa\|_{L^\infty(Q_T)}
+\|b\|_{L^\infty(Q_T)}
\right)
\bigg]
\nonumber\\
&\qquad\quad\times
\left(
\|u_0\|_{H_0^1(0,1)}
+\|f\|_{L^2(0,T;L^2(0,1))}
\right),
\label{second_der}
\end{align}
where \(C>0\) depends only on \(T\) and the interval \((0,1)\).
\end{lem}

\begin{proof}
The basic energy bounds follow from
Theorem~\ref{thm:energy-estimate}. It remains to establish the
higher time--space regularity.

Let \(u_m\in H^1(0,T;V_m)\) be the Galerkin approximation
introduced above. Since
\(\partial_tu_m(t)\in V_m\) for almost every \(t\in(0,T)\), it may
be used as a test function in the Galerkin identity. Taking real
parts gives
\begin{align*}
\|\partial_tu_m\|_{L^2(0,1)}^2
&+\operatorname{Re}
  (a\partial_xu_m,\partial_x\partial_tu_m)\\
&=-\operatorname{Re}
  (b\partial_xu_m,\partial_tu_m)
  -\operatorname{Re}(qu_m,\partial_tu_m)
  +\operatorname{Re}(f,\partial_tu_m).
\end{align*}

Define the weighted energy
\[
E_m(t):=
\int_0^1
a(t,x)|\partial_xu_m(t,x)|^2\,dx.
\]
Then
\[
\operatorname{Re}
(a\partial_xu_m,\partial_x\partial_tu_m)
=
\frac12E_m'(t)
-\frac12\int_0^1
\partial_ta(t,x)|\partial_xu_m(t,x)|^2\,dx.
\]
Therefore,
\begin{align*}
\|\partial_tu_m\|_{L^2(0,1)}^2
+\frac12E_m'(t)
&=
\frac12\int_0^1
\partial_ta\,|\partial_xu_m|^2\,dx\\
&\quad
-\operatorname{Re}
(b\partial_xu_m,\partial_tu_m)
-\operatorname{Re}(qu_m,\partial_tu_m)
+\operatorname{Re}(f,\partial_tu_m).
\end{align*}

By the Cauchy--Schwarz and Young inequalities,
\begin{align*}
\big|
(b\partial_xu_m,\partial_tu_m)
\big|
&\leq
\frac14\|\partial_tu_m\|_{L^2(0,1)}^2
+\|b\|_{L^\infty(Q_T)}^2
 \|\partial_xu_m\|_{L^2(0,1)}^2,\\
\big|
(qu_m,\partial_tu_m)
\big|
&\leq
\frac14\|\partial_tu_m\|_{L^2(0,1)}^2
+\|q\|_{L^\infty(0,1)}^2
 \|u_m\|_{L^2(0,1)}^2,\\
\big|
(f,\partial_tu_m)
\big|
&\leq
\frac14\|\partial_tu_m\|_{L^2(0,1)}^2
+\|f\|_{L^2(0,1)}^2.
\end{align*}
Consequently,
\begin{align}
\frac14\|\partial_tu_m\|_{L^2(0,1)}^2
+\frac12E_m'(t)
&\leq
\left(
\frac12
\|\partial_ta\|_{L^\infty(Q_T)}
+
\|b\|_{L^\infty(Q_T)}^2
\right)
\|\partial_xu_m\|_{L^2(0,1)}^2
\nonumber\\
&\quad
+\|q\|_{L^\infty(0,1)}^2
 \|u_m\|_{L^2(0,1)}^2
+\|f\|_{L^2(0,1)}^2.
\label{eq:weighted-energy}
\end{align}

Integrating \eqref{eq:weighted-energy} over \((0,t)\), taking the
supremum over \(t\in[0,T]\), and also taking \(t=T\), we obtain
\begin{align}
&\|\partial_tu_m\|_{L^2(0,T;L^2(0,1))}^2
+\sup_{t\in[0,T]}E_m(t)
\nonumber\\
&\qquad\leq
C\Big[
E_m(0)
+
\left(
\|\partial_ta\|_{L^\infty(Q_T)}
+\|b\|_{L^\infty(Q_T)}^2
\right)
\|\partial_xu_m\|_{L^2(0,T;L^2(0,1))}^2
\nonumber\\
&\hspace{3cm}
+\|q\|_{L^\infty(0,1)}^2
 \|u_m\|_{L^2(0,T;L^2(0,1))}^2
+\|f\|_{L^2(0,T;L^2(0,1))}^2
\Big].
\label{eq:weighted-integrated}
\end{align}

Since \(a(t,x)\geq\alpha\),
\[
E_m(t)
\geq
\alpha\|\partial_xu_m(t)\|_{L^2(0,1)}^2.
\]
Since \(a\in W^{1,\infty}(Q_T)\), it admits a Lipschitz
representative on \(\overline{Q_T}\), and hence \(a(0,\cdot)\)
is well defined. Moreover, \(u_m(0)=P_mu_0\). The stability of
the sine projection \(P_m\) in \(H_0^1(0,1)\) yields
\[
E_m(0)
\leq
\|a(0,\cdot)\|_{L^\infty(0,1)}
\|\partial_xP_mu_0\|_{L^2(0,1)}^2
\leq
\|a\|_{L^\infty(Q_T)}
\|u_0\|_{H_0^1(0,1)}^2.
\]

Combining \eqref{eq:weighted-integrated} with the basic energy
estimate gives a bound, uniform in \(m\), for
\[
\partial_tu_m
\quad\text{in }L^2(0,T;L^2(0,1))
\]
and for
\[
\partial_xu_m
\quad\text{in }L^\infty(0,T;L^2(0,1)).
\]
By the basic energy estimate and \eqref{eq:weighted-integrated},
the sequences \(\{u_m\}\), \(\{\partial_tu_m\}\), and
\(\{\partial_xu_m\}\) are bounded in the corresponding spaces.
Hence, after passing to a subsequence,
\[
u_m\rightharpoonup u
\quad\text{in }L^2(0,T;H_0^1(0,1)),
\]
\[
\partial_tu_m\rightharpoonup w
\quad\text{in }L^2(0,T;L^2(0,1)),
\]
and
\[
\partial_xu_m\stackrel{*}{\rightharpoonup}z
\quad\text{in }L^\infty(0,T;L^2(0,1)).
\]
Here \(\rightharpoonup\) denotes weak convergence, whereas
\(\stackrel{*}{\rightharpoonup}\) denotes weak-star convergence. The limit \(u\) is the unique weak solution given by
Theorem~\ref{thm:1D-wellposedness}. The standard distributional
identification of weak derivatives gives
\[
w=\partial_tu,
\qquad
z=\partial_xu.
\]
Therefore,
\[
\partial_tu\in L^2(0,T;L^2(0,1)),
\qquad
\partial_xu\in L^\infty(0,T;L^2(0,1)).
\]
By weak and weak-\(*\) lower semicontinuity, the corresponding
bounds pass to the limit.

Using the equation,
\[
\partial_x(a\partial_xu)
=
\partial_tu+b\partial_xu+qu-f,
\]
we obtain
\begin{align*}
\|\partial_x(a\partial_xu)\|_{L^2(0,T;L^2(0,1))}
&\leq
\|\partial_tu\|_{L^2(0,T;L^2(0,1))}
+\|b\|_{L^\infty(Q_T)}
 \|\partial_xu\|_{L^2(0,T;L^2(0,1))}\\
&\quad
+\|q\|_{L^\infty(0,1)}
 \|u\|_{L^2(0,T;L^2(0,1))}
+\|f\|_{L^2(0,T;L^2(0,1))}.
\end{align*}
Thus,
\[
\partial_x(a\partial_xu)
\in L^2(0,T;L^2(0,1)).
\]

Finally,
\[
a\,\partial_x^2u
=
\partial_x(a\partial_xu)
-(\partial_xa)\partial_xu.
\]
Since \(a\geq\alpha\), it follows that
\begin{align*}
\alpha
\|\partial_x^2u\|_{L^2(0,T;L^2(0,1))}
\leq
\|\partial_x(a\partial_xu)\|_{L^2(0,T;L^2(0,1))}
+\|\partial_xa\|_{L^\infty(Q_T)}
 \|\partial_xu\|_{L^2(0,T;L^2(0,1))}.
\end{align*}
Collecting these estimates proves \eqref{second_der}.
\end{proof}

The next section applies these coefficient-dependent estimates to
regularising nets. The classical variational theory recalled above is
used only as an input and is not part of the paper's novelty.

\section{Very weak solutions with nonhomogeneous Dirichlet data}
\label{sec:very-weak}

Let
\[
Q_T:=(0,T)\times(0,1).
\]
We consider the parabolic problem
\begin{equation}
\label{5.2}
\left\{
\begin{array}{l}
\partial_tu-\partial_x(a(t,x)\partial_xu)+b(t,x)\partial_xu+q(x)u
    =f(t,x), \qquad (t,x)\in Q_T,\\
u(0,x)=u_0(x), \qquad x\in(0,1),\\
u(t,0)=g_0(t),\quad
u(t,1)=g_1(t), \qquad t\in(0,T).
\end{array}
\right.
\end{equation}
The drift coefficient, potential, source, initial datum, and boundary
data may be distributions. The diffusion coefficient is allowed to be
singular within the uniformly positive class
\[
a=\alpha+\mu,
\qquad
\alpha>0,
\qquad
\mu\in\mathcal M_+(Q_T),
\qquad
\mu(Q_T)<\infty,
\]
where \(\mathcal M_+(Q_T)\) denotes the space of nonnegative Radon
measures on \(Q_T\). This restriction guarantees that the regularised
problems remain uniformly parabolic. Since \(\mu\) is nonnegative, we write
\[
\|\mu\|_{\mathcal M}:=\mu(Q_T).
\]

Throughout this section, the coefficients and their regularisations
are assumed to be real-valued.

For an open set \(\Omega\subset\mathbb R^d\), we use the notation
\[
\mathcal D(\Omega):=C_c^\infty(\Omega),
\qquad
\mathcal D'(\Omega):=(\mathcal D(\Omega))'.
\]

The notions of moderate and negligible nets originate in the theory of
Colombeau-type generalised functions; see \cite{Ober}. The
regularisation-based notion of a very weak solution used below follows
the framework developed for partial differential equations with
irregular coefficients in \cite{GR15}. We recall the required
definitions and adapt them to the present parabolic problem.

\begin{defn}[Moderate and negligible nets]
\label{def:moderate-nets-1D}
Let \(X(\Omega)\) be a normed function space continuously embedded in
\(\mathcal D'(\Omega)\).
\begin{enumerate}
    \item [a)] For \(v\in\mathcal D'(\Omega)\), a net
\((v_\varepsilon)_{\varepsilon\in(0,1]}\subset X(\Omega)\) is called
an \emph{\(X\)-regularisation} of \(v\) if
\[
v_\varepsilon\longrightarrow v
\quad\text{in }\mathcal D'(\Omega)
\quad\text{as }\varepsilon\to0.
\]

\item [b)] A net \(
(v_\varepsilon)_{\varepsilon\in(0,1]}
\subset X(\Omega)
\) is called \emph{\(X\)-moderate} if there exist \(C>0\) and
\(N\in\mathbb N\) such that
\[
\|v_\varepsilon\|_{X(\Omega)}
\leq C\varepsilon^{-N},
\qquad
\varepsilon\in(0,1].
\]

\item [c)] It is called \emph{log-\(X\)-moderate} if there exists \(C>0\)
such that
\[
\|v_\varepsilon\|_{X(\Omega)}
\leq C\bigl(1+|\log\varepsilon|\bigr),
\qquad
\varepsilon\in(0,1].
\]

\item [d)] It is called \emph{\(X\)-negligible} if, for every \(m\in\mathbb N\),
there exists \(C_m>0\) such that
\[
\|v_\varepsilon\|_{X(\Omega)}
\leq C_m\varepsilon^m,
\qquad
\varepsilon\in(0,1].
\]
\end{enumerate}
\end{defn}

\begin{defn}[Uniformly positive distributional diffusivity]
\label{def:positive-diffusivity}
Let \(\alpha>0\). A distribution \(a\in\mathcal D'(Q_T)\) is called
\emph{uniformly positive with lower bound \(\alpha\)} if
\begin{equation}
\label{distributional-positivity}
\langle a-\alpha,\varphi\rangle\geq0
\qquad
\text{for every }
\varphi\in C_c^\infty(Q_T),\quad \varphi\geq0,
\end{equation}
where \(\alpha\) is identified with the corresponding regular
distribution.

By the characterisation of positive distributions
\cite{Friedlander}, condition
\eqref{distributional-positivity} is equivalent to
\begin{equation}
\label{positive-measure-decomposition}
a=\alpha+\mu
\quad\text{in }\mathcal D'(Q_T),
\qquad
\mu\in\mathcal M_+(Q_T).
\end{equation}

Assume additionally that \(\mu(Q_T)<\infty\). A net
\((a_\varepsilon)_{\varepsilon\in(0,1]}\) is called an
\emph{admissible regularisation} of \(a\) if
\[
a_\varepsilon\in W^{1,\infty}(Q_T),
\qquad
a_\varepsilon\longrightarrow a
\quad\text{in }\mathcal D'(Q_T),
\]
the net is \(W^{1,\infty}(Q_T)\)-moderate, and
\begin{equation}
\label{uniform-ellipticity-net}
a_\varepsilon(t,x)\geq\alpha
\qquad
\text{for almost every }(t,x)\in Q_T,
\quad
\varepsilon\in(0,1].
\end{equation}
\end{defn}

\begin{rem}[Positivity-preserving regularisation]
\label{rem:positive-regularisation}
Extend \(\mu\) by zero to \(\mathbb R^2\), and let
\[
\rho\in C_c^\infty(\mathbb R^2),
\qquad
\rho\geq0,
\qquad
\int_{\mathbb R^2}\rho=1.
\]
Define
\[
\rho_\varepsilon(z)
=
\varepsilon^{-2}\rho(z/\varepsilon),
\qquad
a_\varepsilon
=
\alpha+\mu*\rho_\varepsilon.
\]
Then
\[
a_\varepsilon\geq\alpha,
\qquad
a_\varepsilon\longrightarrow a
\quad\text{in }\mathcal D'(Q_T),
\]
and
\[
\|a_\varepsilon\|_{L^\infty(Q_T)}
\leq
\alpha+C\|\mu\|_{\mathcal M}\varepsilon^{-2},
\]
\[
\|\nabla_{t,x}a_\varepsilon\|_{L^\infty(Q_T)}
\leq
C\|\mu\|_{\mathcal M}\varepsilon^{-3}.
\]
Thus, \((a_\varepsilon)\) is an admissible regularisation.

Conversely, if
\[
a_\varepsilon\longrightarrow a
\quad\text{in }\mathcal D'(Q_T),
\qquad
a_\varepsilon\geq\alpha,
\]
then \(a-\alpha\) is a positive distribution and therefore a
nonnegative Radon measure. Hence an arbitrary distribution cannot
admit a uniformly elliptic regularising net.
\end{rem}

Accordingly, the regularisations below are required to belong to
the function spaces explicitly stated in Definition~\ref{defn5.2}; they
are not assumed to be \(C^\infty\) unless this is stated separately.

\begin{defn}[Very weak solution]
\label{defn5.2}
Assume that
\[
a=\alpha+\mu,
\qquad
\alpha>0,
\qquad
\mu\in\mathcal M_+(Q_T),
\qquad
\mu(Q_T)<\infty,
\]
and
\[
b,f\in\mathcal D'(Q_T),
\qquad
q,u_0\in\mathcal D'(0,1),
\qquad
g_0,g_1\in\mathcal D'(0,T).
\]

A net \((u_\varepsilon)_{\varepsilon\in(0,1]}
\) is called a \emph{very weak solution} of \eqref{5.2} if there exist:

\begin{itemize}
\item an admissible regularisation \((a_\varepsilon)\) of \(a\);

\item an
\(L^\infty(0,T;W^{1,\infty}(0,1))\)-moderate regularisation
\((b_\varepsilon)\) of \(b\), satisfying
\[
\|\partial_xb_\varepsilon\|_{L^\infty(Q_T)}
\leq C\bigl(1+|\log\varepsilon|\bigr);
\]

\item an \(L^\infty(0,1)\)-moderate regularisation
\((q_\varepsilon)\) of \(q\), satisfying
\[
\|q_{-,\varepsilon}\|_{L^\infty(0,1)}
\leq C\bigl(1+|\log\varepsilon|\bigr),
\qquad
q_{-,\varepsilon}:=\max\{-q_\varepsilon,0\};
\]

\item moderate regularisations
\[
(f_\varepsilon)\subset L^2(Q_T),
\qquad
(u_{0,\varepsilon})\subset H^1(0,1),
\qquad
(g_{j,\varepsilon})\subset H^1(0,T),
\quad j=0,1;
\]

\item the compatibility conditions
\[
u_{0,\varepsilon}(0)=g_{0,\varepsilon}(0),
\qquad
u_{0,\varepsilon}(1)=g_{1,\varepsilon}(0).
\]
\end{itemize}

For every \(\varepsilon\in(0,1]\), define
\[
\psi_\varepsilon(t,x)
:=
(1-x)g_{0,\varepsilon}(t)
+xg_{1,\varepsilon}(t).
\]
Then \(u_\varepsilon\) is required to satisfy
\[
u_\varepsilon\in L^2(0,T;H^1(0,1)),
\qquad
u_\varepsilon-\psi_\varepsilon
\in L^2(0,T;H_0^1(0,1)),
\]
and
\[
\partial_tu_\varepsilon
\in L^2(0,T;H^{-1}(0,1)),
\]
and to be the unique weak solution of
\begin{equation}
\label{5.3}
\left\{
\begin{array}{l}
\partial_tu_\varepsilon
-\partial_x(a_\varepsilon(t,x)\partial_xu_\varepsilon)
+b_\varepsilon(t,x)\partial_xu_\varepsilon
+q_\varepsilon(x) u_\varepsilon
=f_\varepsilon(t,x), \quad (t,x)\in Q_T,\\
u_\varepsilon(0,x)=u_{0,\varepsilon}(x), \qquad x\in (0,1),\\
u_\varepsilon(t,0)=g_{0,\varepsilon}(t),
\quad
u_\varepsilon(t,1)=g_{1,\varepsilon}(t),
\qquad t\in (0,T).
\end{array}
\right.
\end{equation}
Moreover, the net \((u_\varepsilon)\) is moderate in
\[
H^1(0,T;L^2(0,1))
\cap
L^2(0,T;H^2(0,1)),
\]
equipped with the sum norm.
\end{defn}

\subsection{Boundary lifting and regularised estimates}

With \(\psi_\varepsilon\) as in Definition~\ref{defn5.2}, set
\[
w_\varepsilon:=u_\varepsilon-\psi_\varepsilon.
\]
Then \(w_\varepsilon\) satisfies homogeneous Dirichlet boundary
conditions and solves
\begin{equation}
\label{hom.w}
\left\{
\begin{array}{l}
\partial_tw_\varepsilon
-\partial_x(a_\varepsilon(t,x)\partial_xw_\varepsilon)
+b_\varepsilon(t,x)\partial_xw_\varepsilon
+q_\varepsilon(x)w_\varepsilon
=\widetilde F_\varepsilon(t,x), \quad (t,x)\in Q_T,\\
w_\varepsilon(0,x)
=
u_{0,\varepsilon}(x)
-\psi_\varepsilon(0,x), \qquad x\in (0,1),\\
w_\varepsilon(t,0)
=
w_\varepsilon(t,1)=0, \qquad t\in (0,T),
\end{array}
\right.
\end{equation}
where
\begin{align}
\widetilde F_\varepsilon
=&
f_\varepsilon
-(1-x)g_{0,\varepsilon}'
-xg_{1,\varepsilon}'\notag\\
&+
(\partial_xa_\varepsilon-b_\varepsilon)
(g_{1,\varepsilon}-g_{0,\varepsilon})-q_\varepsilon
\big((1-x)g_{0,\varepsilon}+xg_{1,\varepsilon}\big).
\label{lifted-source}
\end{align}

The embedding
\[
H^1(0,T)\hookrightarrow C([0,T])
\]
gives
\begin{align}
\|\widetilde F_\varepsilon\|_{L^2(Q_T)}
\leq{}&
C\Big[
\|f_\varepsilon\|_{L^2(Q_T)}
+
\big(
1+\|\partial_xa_\varepsilon\|_{L^\infty(Q_T)}
+\|b_\varepsilon\|_{L^\infty(Q_T)}
+\|q_\varepsilon\|_{L^\infty(0,1)}
\big)                                                     \notag\\
&\hspace{35mm}\times
\sum_{j=0}^1
\|g_{j,\varepsilon}\|_{H^1(0,T)}
\Big].
\label{tilde_f_eps}
\end{align}
The compatibility conditions imply
\[
u_{0,\varepsilon}
-\psi_\varepsilon(0,\cdot)
\in H_0^1(0,1),
\]
and
\begin{equation}
\label{lifted-initial-estimate}
\|u_{0,\varepsilon}
-\psi_\varepsilon(0,\cdot)\|_{H_0^1(0,1)}
\leq
C\left(
\|u_{0,\varepsilon}\|_{H^1(0,1)}
+
\sum_{j=0}^1
\|g_{j,\varepsilon}\|_{H^1(0,T)}
\right).
\end{equation}

\begin{prop}[Higher regularity for the regularised problem]
\label{u_eps_high_reg}
Let \(u_\varepsilon\) solve the regularised nonhomogeneous Dirichlet
problem \eqref{5.3}. Assume that
\[
f_\varepsilon\in L^2(Q_T),
\qquad
u_{0,\varepsilon}\in H^1(0,1),
\qquad
g_{0,\varepsilon},g_{1,\varepsilon}\in H^1(0,T),
\]
and that
\[
u_{0,\varepsilon}(0)=g_{0,\varepsilon}(0),
\qquad
u_{0,\varepsilon}(1)=g_{1,\varepsilon}(0).
\]
Suppose further that
\[
a_\varepsilon,\partial_ta_\varepsilon,
\partial_xa_\varepsilon,b_\varepsilon,\partial_xb_\varepsilon
\in L^\infty(Q_T),
\qquad
a_\varepsilon\geq\alpha>0,
\]
and
\[
q_\varepsilon\in L^\infty(0,1),
\qquad
q_{-,\varepsilon}:=\max\{-q_\varepsilon,0\}.
\]

Then
\[
u_\varepsilon
\in H^1(0,T;L^2(0,1))
\cap L^2(0,T;H^2(0,1)),
\]
and
\[
\partial_x
\big(a_\varepsilon\partial_xu_\varepsilon\big)
\in L^2(Q_T).
\]

Set
\[
B_\varepsilon
:=
\|\partial_xb_\varepsilon\|_{L^\infty(Q_T)}
+
2\|q_{-,\varepsilon}\|_{L^\infty(0,1)}
\]
and
\begin{align*}
K_\varepsilon
:=
&1
+\|a_\varepsilon\|_{L^\infty(Q_T)}^{1/2}
+\|q_\varepsilon\|_{L^\infty(0,1)}
\\&+
\frac{1}{\sqrt{\alpha}}
\left(
\|\partial_ta_\varepsilon\|_{L^\infty(Q_T)}^{1/2}
+\|\partial_xa_\varepsilon\|_{L^\infty(Q_T)}
+\|b_\varepsilon\|_{L^\infty(Q_T)}
\right).
\end{align*}
Then
\begin{align}
&
\|\partial_tu_\varepsilon\|_{L^2(Q_T)}
+
\sqrt{\alpha}\,
\|\partial_xu_\varepsilon\|_{L^\infty(0,T;L^2(0,1))}
+
\left\|
\partial_x
\big(a_\varepsilon\partial_xu_\varepsilon\big)
\right\|_{L^2(Q_T)}
\notag\\
&\qquad
+
\alpha\|\partial_x^2u_\varepsilon\|_{L^2(Q_T)}
\notag\\
&\quad\leq
C e^{TB_\varepsilon}
K_\varepsilon
\left[
\|u_{0,\varepsilon}
-\psi_\varepsilon(0,\cdot)\|_{H_0^1(0,1)}
+
\|\widetilde F_\varepsilon\|_{L^2(Q_T)}
\right]
\notag\\
&\qquad+
C\left(
1+\sqrt{\alpha}
+\|\partial_xa_\varepsilon\|_{L^\infty(Q_T)}
\right)
\sum_{j=0}^{1}
\|g_{j,\varepsilon}\|_{H^1(0,T)},
\label{D_t_u_L2}
\end{align}
where \(C>0\) depends only on \(T\) and the spatial interval
\((0,1)\).
\end{prop}

\begin{proof}
The function
\[
w_\varepsilon
=
u_\varepsilon-\psi_\varepsilon
\]
solves the homogeneous problem \eqref{hom.w}. The compatibility
conditions give
\[
w_\varepsilon(0,\cdot)
=
u_{0,\varepsilon}
-\psi_\varepsilon(0,\cdot)
\in H_0^1(0,1).
\]
Applying Lemma~\ref{lem:time-space-second-derivative} to
\(w_\varepsilon\) yields
\begin{align*}
&
\|\partial_tw_\varepsilon\|_{L^2(Q_T)}
+
\sqrt{\alpha}\,
\|\partial_xw_\varepsilon\|_{L^\infty(0,T;L^2(0,1))}
+
\left\|
\partial_x
\big(a_\varepsilon\partial_xw_\varepsilon\big)
\right\|_{L^2(Q_T)}
\\
&\qquad
+
\alpha\|\partial_x^2w_\varepsilon\|_{L^2(Q_T)}
\\
&\quad\leq
C e^{TB_\varepsilon}
K_\varepsilon
\left[
\|w_\varepsilon(0,\cdot)\|_{H_0^1(0,1)}
+
\|\widetilde F_\varepsilon\|_{L^2(Q_T)}
\right].
\end{align*}
Moreover,
\[
\sqrt{\alpha}\,
\|\partial_x\psi_\varepsilon\|_{L^\infty(0,T;L^2(0,1))}
\leq
C\sqrt{\alpha}
\sum_{j=0}^{1}
\|g_{j,\varepsilon}\|_{H^1(0,T)}.
\]

Since
\[
u_\varepsilon=w_\varepsilon+\psi_\varepsilon,
\qquad
\partial_x^2\psi_\varepsilon=0,
\]
we have
\[
\partial_tu_\varepsilon
=
\partial_tw_\varepsilon+\partial_t\psi_\varepsilon
\]
and
\[
\partial_x
\big(a_\varepsilon\partial_xu_\varepsilon\big)
=
\partial_x
\big(a_\varepsilon\partial_xw_\varepsilon\big)
+
(\partial_xa_\varepsilon)
(g_{1,\varepsilon}-g_{0,\varepsilon}).
\]
Moreover,
\[
\|\partial_t\psi_\varepsilon\|_{L^2(Q_T)}
\leq
C\sum_{j=0}^{1}
\|g_{j,\varepsilon}\|_{H^1(0,T)}
\]
and
\[
\left\|
(\partial_xa_\varepsilon)
(g_{1,\varepsilon}-g_{0,\varepsilon})
\right\|_{L^2(Q_T)}
\leq
C\|\partial_xa_\varepsilon\|_{L^\infty(Q_T)}
\sum_{j=0}^{1}
\|g_{j,\varepsilon}\|_{H^1(0,T)}.
\]
Combining these inequalities proves \eqref{D_t_u_L2}.
\end{proof}
\begin{thm}[Existence of very weak solutions]
\label{thm5.1}
Assume that the coefficients and data satisfy
\[
a=\alpha+\mu,\qquad
\alpha>0,\qquad
\mu\in\mathcal M_+(Q_T),\qquad
\mu(Q_T)<\infty,
\]
and
\[
b,f\in\mathcal D'(Q_T),\qquad
q,u_0\in\mathcal D'(0,1),\qquad
g_0,g_1\in\mathcal D'(0,T).
\]
Suppose that these distributions admit regularising families
satisfying all the moderateness, logarithmic-growth, compatibility,
and uniform ellipticity conditions imposed on the coefficient and
data regularisations in Definition~\ref{defn5.2}.
Then problem \eqref{5.2} admits a very weak solution.
\end{thm}

\begin{proof}
For every fixed \(\varepsilon\in(0,1]\), the lifting
\(w_\varepsilon=u_\varepsilon-\psi_\varepsilon\) reduces
\eqref{5.3} to the homogeneous problem \eqref{hom.w}. Therefore,
the classical theory recalled in Section~2 gives a unique weak
solution \(u_\varepsilon\).

By the logarithmic assumptions,
\[
B_\varepsilon
\leq C\bigl(1+|\log\varepsilon|\bigr).
\]
Consequently,
\[
e^{TB_\varepsilon}
\leq C\varepsilon^{-N}
\]
for some \(N\in\mathbb N\). By the assumptions on the regularising families, \eqref{tilde_f_eps}, and \eqref{lifted-initial-estimate}, every factor on the right-hand side of \eqref{D_t_u_L2} is moderate. Therefore, Proposition~\ref{u_eps_high_reg} yields
\[
\|\partial_tu_\varepsilon\|_{L^2(Q_T)}
+
\|\partial_xu_\varepsilon\|_{L^\infty(0,T;L^2(0,1))}
+
\|\partial_x^2u_\varepsilon\|_{L^2(Q_T)}
\leq C\varepsilon^{-N'}
\]
for suitable \(C>0\) and \(N'\in\mathbb N\). 
Furthermore, since
\[
w_\varepsilon
=
u_\varepsilon-\psi_\varepsilon
\in L^2(0,T;H_0^1(0,1)),
\]
Poincaré's inequality and the estimates for the lifting imply
\[
\|u_\varepsilon\|_{L^2(Q_T)}
+
\|\partial_xu_\varepsilon\|_{L^2(Q_T)}
\leq
C\left(
\|\partial_xw_\varepsilon\|_{L^2(Q_T)}
+
\sum_{j=0}^{1}
\|g_{j,\varepsilon}\|_{H^1(0,T)}
\right).
\]
It follows that
\[
\|u_\varepsilon\|_{H^1(0,T;L^2(0,1))}
+
\|u_\varepsilon\|_{L^2(0,T;H^2(0,1))}
\leq C\varepsilon^{-N''}
\]
for some \(N''\in\mathbb N\). Therefore
\((u_\varepsilon)\) is a very weak solution.
\end{proof}

\begin{thm}[Uniqueness with respect to equivalent regularisations]
\label{thm5.2}
Let \((u_\varepsilon)\) and
\((\widetilde u_\varepsilon)\) be two very weak solutions generated
by regularising families satisfying the coefficient and data
conditions of Definition~\ref{defn5.2}, with the same ellipticity
lower bound \(\alpha\).

Assume that
\[
a_\varepsilon-\widetilde a_\varepsilon
\quad\text{is }
W^{1,\infty}(Q_T)\text{-negligible},
\]
\[
b_\varepsilon-\widetilde b_\varepsilon
\quad\text{is }
L^\infty(0,T;W^{1,\infty}(0,1))\text{-negligible},
\]
\[
q_\varepsilon-\widetilde q_\varepsilon
\quad\text{is }
L^\infty(0,1)\text{-negligible},
\]
and that the differences of the source, initial datum, and boundary
data are negligible in
\[
L^2(Q_T),\qquad H^1(0,1),\qquad H^1(0,T),
\]
respectively. Then
\[
u_\varepsilon-\widetilde u_\varepsilon
\]
is negligible in
\[
H^1(0,T;L^2(0,1))
\cap
L^2(0,T;H^2(0,1)).
\]
\end{thm}

\begin{proof}
Set
\[
U_\varepsilon
:=
u_\varepsilon-\widetilde u_\varepsilon.
\]
In applying Proposition~\ref{u_eps_high_reg} to the difference
problem, \(B_\varepsilon\) and \(K_\varepsilon\) are defined using
the coefficients \(a_\varepsilon,b_\varepsilon,q_\varepsilon\)
placed on the left-hand side.

Using the coefficients \(a_\varepsilon,b_\varepsilon,q_\varepsilon\)
on the left-hand side, the source term in the equation for
\(U_\varepsilon\) is
\begin{align*}
\Phi_\varepsilon
={}&
-\partial_x\left(
(\widetilde a_\varepsilon-a_\varepsilon)
\partial_x\widetilde u_\varepsilon
\right)\\
&+
(\widetilde b_\varepsilon-b_\varepsilon)
\partial_x\widetilde u_\varepsilon
+
(\widetilde q_\varepsilon-q_\varepsilon)
\widetilde u_\varepsilon
+
f_\varepsilon-\widetilde f_\varepsilon.
\end{align*}
Moreover,
\begin{align*}
\|\Phi_\varepsilon\|_{L^2(Q_T)}
\leq{}&
C\|\widetilde a_\varepsilon-a_\varepsilon\|_{W^{1,\infty}(Q_T)}
\|\widetilde u_\varepsilon\|_{L^2(0,T;H^2(0,1))}\\
&+
\|\widetilde b_\varepsilon-b_\varepsilon\|_{L^\infty(Q_T)}
\|\partial_x\widetilde u_\varepsilon\|_{L^2(Q_T)}\\
&+
\|\widetilde q_\varepsilon-q_\varepsilon\|_{L^\infty(0,1)}
\|\widetilde u_\varepsilon\|_{L^2(Q_T)}
+
\|f_\varepsilon-\widetilde f_\varepsilon\|_{L^2(Q_T)}.
\end{align*}

The solution net \((\widetilde u_\varepsilon)\) is moderate, whereas
all coefficient and data differences are negligible. Since the
product of a negligible scalar net and a moderate net is negligible,
it follows that
\[
(\Phi_\varepsilon)_{\varepsilon\in(0,1]}
\quad\text{is \(L^2(Q_T)\)-negligible}.
\]

Let
\[
G_{j,\varepsilon}
:=
g_{j,\varepsilon}-\widetilde g_{j,\varepsilon},
\qquad j=0,1,
\]
and define the corresponding lifting
\[
\Psi_\varepsilon(t,x)
:=
(1-x)G_{0,\varepsilon}(t)
+xG_{1,\varepsilon}(t).
\]
Applying \eqref{tilde_f_eps} to the difference problem, with \(\Phi_\varepsilon\) in place of \(f_\varepsilon\), and using the moderateness of the coefficient nets together with the \(H^1(0,T)\)-negligibility of \(G_{0,\varepsilon}\) and \(G_{1,\varepsilon}\), we conclude that the source obtained after homogenising the boundary conditions is
\(L^2(Q_T)\)-negligible. Moreover,
\eqref{lifted-initial-estimate} shows that
\[
U_\varepsilon(0,\cdot)-\Psi_\varepsilon(0,\cdot)
\]
is \(H_0^1(0,1)\)-negligible. The compatibility conditions follow by
subtracting those satisfied by the two regularising families.

Finally, the factor
\[
e^{TB_\varepsilon}K_\varepsilon
\]
is moderate. 
Applying Proposition~\ref{u_eps_high_reg} shows that
\(\partial_tU_\varepsilon\), \(\partial_xU_\varepsilon\), and
\(\partial_x^2U_\varepsilon\) are negligible in their corresponding
norms. Since
\[
U_\varepsilon-\Psi_\varepsilon
\in L^2(0,T;H_0^1(0,1)),
\]
Poincaré's inequality and the negligibility of
\(\Psi_\varepsilon\) imply that \(U_\varepsilon\) itself is
\(L^2(Q_T)\)-negligible. Therefore \(U_\varepsilon\) is negligible
in the asserted intersection space.
\end{proof}

\begin{thm}[Consistency]
\label{thm5.3}
Assume that
\[
a\in W^{1,\infty}(Q_T),
\qquad
a(t,x)\geq\alpha_0>0
\quad\text{a.e. in }Q_T,
\]
\[
b\in L^\infty(0,T;W^{1,\infty}(0,1)),
\qquad
q\in L^\infty(0,1),
\qquad
f\in L^2(Q_T),
\]
and
\[
u_0\in H^1(0,1),
\qquad
g_0,g_1\in H^1(0,T),
\]
with the compatibility conditions
\[
u_0(0)=g_0(0),
\qquad
u_0(1)=g_1(0).
\]

Let \(u\) be the unique solution of \eqref{5.2} satisfying
\[
u\in H^1(0,T;L^2(0,1))
\cap L^2(0,T;H^2(0,1)).
\] 
Suppose that the regularising
families satisfy, for some \(0<\alpha\leq\alpha_0\),
\[
a_\varepsilon\geq\alpha
\quad\text{uniformly in }\varepsilon,
\]
and
\[
a_\varepsilon\longrightarrow a
\quad\text{in }W^{1,\infty}(Q_T),
\]
\[
b_\varepsilon\longrightarrow b
\quad\text{in }
L^\infty(0,T;W^{1,\infty}(0,1)),
\qquad
q_\varepsilon\longrightarrow q
\quad\text{in }L^\infty(0,1),
\]
\[
f_\varepsilon\longrightarrow f
\quad\text{in }L^2(Q_T),
\qquad
u_{0,\varepsilon}\longrightarrow u_0
\quad\text{in }H^1(0,1),
\]
and
\[
g_{j,\varepsilon}\longrightarrow g_j
\quad\text{in }H^1(0,T),
\qquad j=0,1.
\]

Assume additionally that, for every \(\varepsilon\in(0,1]\),
\[
u_{0,\varepsilon}(0)=g_{0,\varepsilon}(0),
\qquad
u_{0,\varepsilon}(1)=g_{1,\varepsilon}(0).
\]
Then the corresponding regularised solutions satisfy
\[
u_\varepsilon\longrightarrow u
\]
in
\[
H^1(0,T;L^2(0,1))
\cap
L^2(0,T;H^2(0,1)).
\]
\end{thm}

\begin{proof}
Set
\[
U_\varepsilon:=u_\varepsilon-u.
\]
Here \(B_\varepsilon\) and \(K_\varepsilon\) are defined using
\(a_\varepsilon,b_\varepsilon,q_\varepsilon\). Then \(U_\varepsilon\) satisfies a regular parabolic equation whose
source is
\begin{align*}
\Phi_\varepsilon
={}&
\partial_x\left(
(a_\varepsilon-a)\partial_xu
\right)
-
(b_\varepsilon-b)\partial_xu\\
&-
(q_\varepsilon-q)u
+
f_\varepsilon-f.
\end{align*}
Since
\[
u\in L^2(0,T;H^2(0,1))
\]
and the coefficients converge in the stated norms,
\[
\|\Phi_\varepsilon\|_{L^2(Q_T)}
\longrightarrow0.
\]

Let
\[
G_{j,\varepsilon}
:=
g_{j,\varepsilon}-g_j,
\qquad j=0,1,
\]
and let
\[
\Psi_\varepsilon(t,x)
=
(1-x)G_{0,\varepsilon}(t)
+xG_{1,\varepsilon}(t).
\]
Applying \eqref{tilde_f_eps} to the difference problem, with \(\Phi_\varepsilon\) in place of \(f_\varepsilon\), and using the uniform boundedness of the coefficient families together with
\[
G_{j,\varepsilon}\longrightarrow0
\quad\text{in }H^1(0,T),\qquad j=0,1,
\]
we conclude that the source obtained after homogenising the boundary conditions converges to zero in \(L^2(Q_T)\). Likewise,
\eqref{lifted-initial-estimate} gives
\[
\|U_\varepsilon(0,\cdot)
-\Psi_\varepsilon(0,\cdot)\|_{H_0^1(0,1)}
\longrightarrow0.
\]

Furthermore, \(e^{TB_\varepsilon}K_\varepsilon\) is uniformly bounded,
because the coefficient regularisations converge strongly in the
stated norms and retain the common ellipticity lower bound.
Applying Proposition~\ref{u_eps_high_reg} shows that
\[
\|\partial_tU_\varepsilon\|_{L^2(Q_T)}
+
\|\partial_xU_\varepsilon\|_{L^\infty(0,T;L^2(0,1))}
+
\|\partial_x^2U_\varepsilon\|_{L^2(Q_T)}
\longrightarrow0.
\]
Since
\[
U_\varepsilon-\Psi_\varepsilon
\in L^2(0,T;H_0^1(0,1)),
\]
Poincaré's inequality, together with
\(\Psi_\varepsilon\to0\) in
\(H^1(0,T;H^1(0,1))\), gives
\[
U_\varepsilon\longrightarrow0
\]
in
\[
H^1(0,T;L^2(0,1))
\cap L^2(0,T;H^2(0,1)).
\]
\end{proof}

\begin{rem}[Interpretation of the solution net]
A very weak solution is represented by a moderate net of solutions to
regularised problems. Moderateness provides polynomial control as
\(\varepsilon\to0\), but does not by itself imply convergence in a
classical function space or in \(\mathcal D'(Q_T)\).

Theorem~\ref{thm5.3} shows that, for sufficiently regular
coefficients and data, the regularised solutions converge to the
unique regular solution. For genuinely distributional coefficients or
data, however, moderateness alone does not imply convergence or
association with a distribution. Such behaviour must be investigated
separately.
\end{rem}

\section{Numerical Experiments}
In this Section, we carry out numerical experiments for parabolic equations with singular coefficients, mainly in the one-dimensional case. We also present two-dimensional simulations for Cases~2 and~4 in order to illustrate the spatial effects of singular diffusion and potential coefficients.

More precisely, we consider the initial–boundary value problem
\begin{align}\label{n1.1}
    \partial_t u(t,x) - \partial_x \left(a(t,x)\partial_x u(t,x)\right) 
    + b(t,x) \partial_x u(t,x) 
    + q(x) u(t,x) 
    = f(t,x),\nonumber\\
    \quad (t,x) \in [0,T] \times (0,1),
\end{align}
with initial condition
\begin{equation}\label{n1.2}
    u(0,x) = u_0(x), \quad x \in (0,1),
\end{equation}
and inhomogeneous Dirichlet boundary conditions
\begin{equation}\label{n1.3}
    u(t,0) = g_0(t), \quad u(t,1) = g_1(t), \quad t \in [0,T].
\end{equation}

We are interested in singular drift, potential, initial, and boundary data, which may contain $\delta$-type contributions. In accordance with
Definition~\ref{def:positive-diffusivity}, the singular diffusivity is restricted
to
\[
a=\alpha+\mu,\qquad \alpha>0,\qquad \mu\ge0,
\]
so that it remains strictly positive after regularisation. In the numerical
examples below, \(\alpha=1\) and \(\mu\) is either zero or a Dirac measure. As it was theoretically outlined in \cite{RT17a} and \cite{RT17b}, we start to analyse our problem by regularising a distributional valued functions $a(t, x), b(t, x), q(t), u_0(x), g_0(t), g_1(t)$  by a parameter $\varepsilon$, that is, we set
\begin{align}
&a_{\varepsilon}(x):=(a*\varphi_{\varepsilon})(x), \quad b_{\varepsilon}(t):=(b*\varphi_{\varepsilon})(t), \quad q_{\varepsilon}(x):=(q*\varphi_{\varepsilon})(x), \quad \notag \\ 
&g_{0,\varepsilon}(t):=(g_0*\varphi_{\varepsilon})(t),\quad g_{1,\varepsilon}(t):=(g_1*\varphi_{\varepsilon})(t)
\end{align}

as the convolution with the mollifier
\begin{equation}
\varphi_{\varepsilon}(x)=\frac{1}{\varepsilon}\varphi(x/\varepsilon),
\end{equation}
where
$\varphi(x)=
c \exp{\left(\frac{1}{x^{2}-1}\right)}$ for $|x| < 1,$ and $\varphi(x)=0$ otherwise. Here  $c \simeq 2.2523$ to get
$\int\limits_{-\infty}^{\infty}  \varphi(x)dx=1.$
For the time-dependent coefficients and data we use the analogous mollifier in
the time variable. Since \(\varphi_\varepsilon\ge0\) and \(\mu\ge0\),
\[
a_\varepsilon(x)=\alpha+(\mu*\varphi_\varepsilon)(x)\ge\alpha>0,
\]
so every numerical regularisation is uniformly parabolic.
Then, instead of \eqref{1.1} we consider the regularised problem

\begin{align}\label{Equation_num2}
\small
\partial_t u_{\varepsilon}(t,x) - \partial_x \left(a_{\varepsilon}(t,x)\partial_x u_{\varepsilon}(t,x)\right)
& + b_{\varepsilon}(t,x) \partial_x u_{\varepsilon}(t,x)
+ q_{\varepsilon}(x) u_{\varepsilon}(t,x)\\
&= f_{\varepsilon}(t,x),\nonumber\
\, (t,x) \in [0,T] \times (0,1),
\end{align}
with
\[
u_{\varepsilon}(0,x)=u_{0,\varepsilon}(x),
\]
and
\[
u_\varepsilon(t,0)=g_{0,\varepsilon}(t),
\qquad
u_\varepsilon(t,1)=g_{1,\varepsilon}(t).
\]
Here, we put
\[
u_0(x) = 5,\exp\left(-\frac{(x-0.5)^2}{0.025}\right),
\]
and $f=0$ in the numerical examples considered below

\subsection{Numerical discretisation and verification}

We briefly describe the numerical method used to solve the regularised problem. Let
\[
x_i=ih,\qquad i=0,\ldots,N,\qquad h=\frac1N,
\]
and
\[
t_n=n\tau,\qquad n=0,\ldots,M,\qquad \tau=\frac{T}{M}.
\]
We denote by
\[
u_i^n\approx u_\varepsilon(t_n,x_i)
\]
the numerical approximation to the regularised solution.

For every fixed $\varepsilon>0$, all regularised coefficients are smooth. The time derivative is approximated by the backward Euler formula, the second-order spatial derivative by the standard central difference, and the first-order spatial derivatives by forward differences. Writing
\[
-\partial_x(a_\varepsilon u_x) =
-a_\varepsilon u_{xx}
-(\partial_xa_\varepsilon)u_x,
\]
the resulting semi-implicit finite-difference scheme is

\begin{align}\label{num_scheme}
\frac{u_i^{n+1}-u_i^n}{\tau}
&-
a_{\varepsilon,i}
\frac{u_{i+1}^{n+1}-2u_i^{n+1}+u_{i-1}^{n+1}}{h^2}
\nonumber\
-
\frac{a_{\varepsilon,i+1}-a_{\varepsilon,i}}{h}
\frac{u_{i+1}^n-u_i^n}{h} \\
&+ b_{\varepsilon,n}
\frac{u_{i+1}^n-u_i^n}{h}
+
q_{\varepsilon,i}u_i^{n+1}
= f_{\varepsilon,i}^n,
\end{align}

for $i=1,\ldots,N-1$, where
\[
a_{\varepsilon,i}=a_\varepsilon(x_i),\qquad
b_{\varepsilon,n}=b_\varepsilon(t_n),\qquad
q_{\varepsilon,i}=q_\varepsilon(x_i).
\]
The initial and boundary conditions are imposed directly as
\[
u_i^0=u_{0,\varepsilon}(x_i),\qquad
u_0^n=g_{0,\varepsilon}(t_n),\qquad
u_N^n=g_{1,\varepsilon}(t_n).
\]

The above discretisation is a standard semi-implicit finite-difference approach for linear parabolic equations, based on the backward Euler approximation in time and finite differences in space. The stability and convergence properties of such implicit finite-difference discretisations for parabolic problems are classical; see, for example, \cite{MM2005,Strikwerda2004}. In the present scheme, the principal diffusion and reaction terms are treated implicitly, while the lower-order first-order terms are evaluated explicitly. This preserves the computational simplicity of the method while retaining the stability properties associated with the implicit treatment of the principal parabolic part. At every time step, the implicit part of \eqref{num_scheme} leads to a tridiagonal linear system, which is solved by the Thomas algorithm.

The singular coefficients and data are not inserted directly into the numerical scheme. A Dirac distribution concentrated at $z=z_*$, where $z$ denotes either the spatial or temporal variable, is approximated by
\[
\delta_{\varepsilon,z_*}(z)
\frac{1}{\varepsilon}
\varphi\left(\frac{z-z_*}{\varepsilon}\right),
\]
where $\varphi$ is the compactly supported mollifier introduced above. Thus, for example,
\[
a_\varepsilon(x) =
1+\delta_{\varepsilon,0.45}(x)
\]
in Case~2, and the other singular coefficients and boundary data are treated analogously.

For smooth regularised data, the above discretisation is consistent, being first-order accurate in time and, due to the one-sided approximation of the first-order terms, first-order accurate in space. Moreover, since $a_\varepsilon\ge1$ in all the considered examples, the implicit diffusion part preserves the parabolic character of the regularised problem. To ensure that the regularised singularities are properly resolved, the spatial and temporal mesh sizes are chosen sufficiently small compared with $\varepsilon$, namely
\[
h\ll\varepsilon
\quad\text{for spatial singularities},\qquad
\tau\ll\varepsilon
\quad\text{for temporal singularities}.
\]
For each fixed $\varepsilon$, the numerical results are additionally checked under mesh refinement before investigating the dependence of the solution on the regularisation parameter. This allows the effects caused by the singular coefficients to be distinguished from numerical artefacts due to insufficient resolution.

\subsection{One-dimensional experiments}

The purpose of the following examples is not merely to solve the regularised problems for several fixed values of $\varepsilon$, but also to illustrate numerically the behaviour of the regularising net $(u_\varepsilon)_\varepsilon$ appearing in the construction of the very weak solution. For every value of $\varepsilon$, the corresponding singular coefficient or boundary datum is first replaced by its mollified representative and the complete regularised problem is then solved independently using the same numerical procedure described above. The computations are repeated for decreasing values of $\varepsilon$ in order to observe how the concentration of the regularised singularity influences the corresponding solution.

As $\varepsilon$ decreases, the support of $\delta_{\varepsilon,z_*}$ shrinks towards the singular point $z_*$, while its total mass remains equal to one. Consequently, its maximum amplitude increases as its effect becomes concentrated in an increasingly small spatial or temporal neighbourhood. The numerical comparisons below are therefore intended to show whether the influence of the singular coefficient becomes localised near its support and how the corresponding solution profiles behave away from this region. We emphasise that the general notion of a very weak solution is formulated in terms of the regularising net $(u_\varepsilon)_\varepsilon$ and does not, in general, require convergence of this net in a classical function space as $\varepsilon\to0$. Thus, the experiments below illustrate the behaviour of the regularised representatives rather than provide a proof of classical convergence.

\begin{figure}[h]
\begin{minipage}[h]{\linewidth}
\center{\includegraphics[scale=0.3]{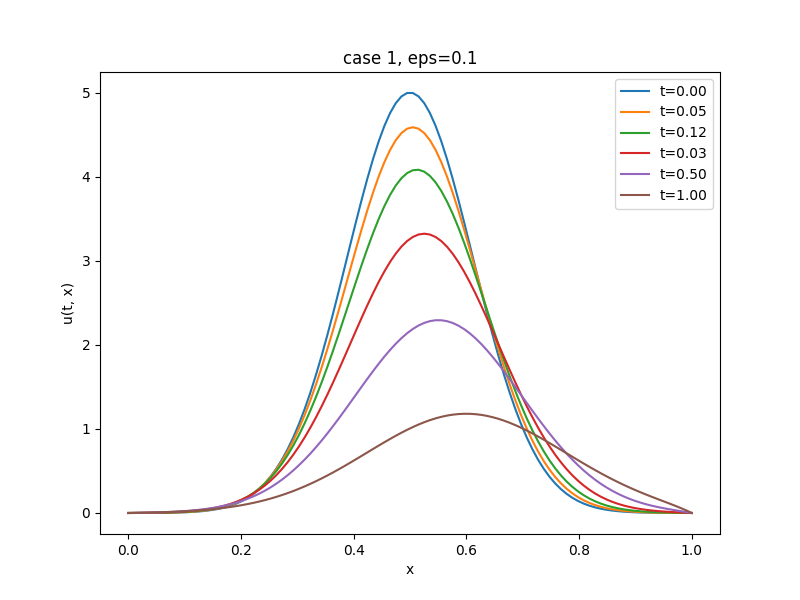}}
\end{minipage}
\caption{Case 1}
\label{case1}
\end{figure}

For the coefficients $a, b, q, g_0, g_1$ we consider the following cases, with $\delta$ denoting the standard Dirac’s delta-distribution:
First, we study the following model situation:
\begin{itemize}
\item Case 1, $a = 1, \quad b = 1,  \quad q = 1, \quad g_0 = g_1 = 0$;
\item Case 2, $a = 1 + \delta(x-0.45), \quad b = 1, \quad q = 1, \quad g_0 = g_1 = 0$;
\item Case 3, $a = 1, \quad b = 1 + \delta(t-0.5), \quad q = 1, \quad g_0 = g_1 = 0$;
\item Case 4, $a = 1, \quad b = 1, \quad q = 1 + \delta(x-0.45),  \quad g_0 = g_1 = 0$;
\item Case 5, $a = 1, \quad b = 1,  \quad q = 1,\quad g_0 = 0 \quad g_1 = \delta(t-0.45).$
\end{itemize}
Thus, the regularised singular quantities used in Cases~2--5 are, respectively,
\[
a_\varepsilon(x)=1+\varphi_\varepsilon(x-0.45),
\]
\[
b_\varepsilon(t)=1+\varphi_\varepsilon(t-0.5),
\]
\[
q_\varepsilon(x)=1+\varphi_\varepsilon(x-0.45),
\]
and
\[
g_{1,\varepsilon}(t)=\varphi_\varepsilon(t-0.45).
\]
All the remaining coefficients and data are kept unchanged when $\varepsilon$ varies. Therefore, each curve corresponding to a different value of $\varepsilon$ is obtained by solving a separate regularised initial--boundary value problem.
For comparison with other cases in figure \ref{case1}, we illustrate the regular case when all the coefficients are constant, such that $a=1.0, b = 1.0, q = 1, g_0 = g_1 = 0$, the behaviour of the numerical solution $u(t,x)$ at different time steps.
The evolution of the heat distribution in Figure \ref{case1} displays smooth temporal decay accompanied by a rightward shift, which is a direct consequence of the transport effected by the $u_x$ term.

\begin{figure}[h]
\includegraphics[scale=0.3]{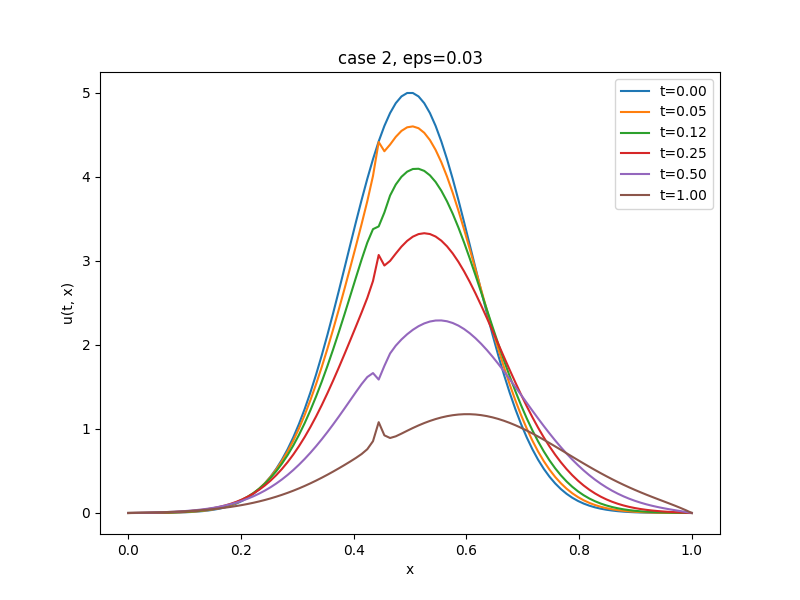}
\includegraphics[scale=0.3]{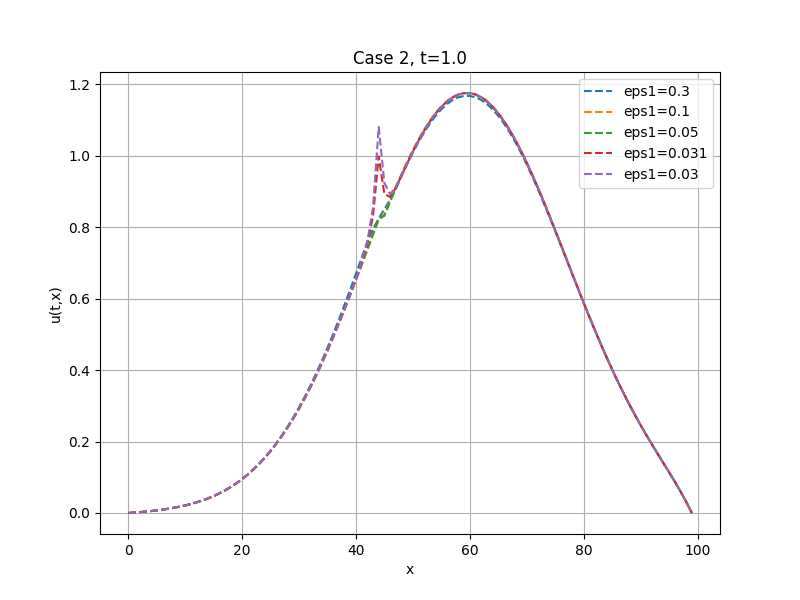}
\caption{Case 2}
\label{case2}
\end{figure}

Figure \ref{case2} illustrates the behaviour of the numerical solution $u(t,x)$, which serves as the solution to equation \ref{Equation_num2}, analysed for $\varepsilon$=0.3, 0.1, 0.05, 0.031 and 0.003 in Case 2. Here the diffusivity is the strictly positive measure
\(a=1+\delta_{0.45}\), not the Dirac measure alone; correspondingly,
\(a_\varepsilon=1+\varphi_\varepsilon(\cdot-0.45)\ge1\).

Figure~\ref{case2} illustrates the temporal decay of the heat distribution $u(t,x)$. The singular diffusion coefficient produces a localised modification of the solution near $x=0.45$. As $\varepsilon$ decreases, the regularised diffusivity becomes increasingly concentrated around this point, and the region in which the corresponding solution is most strongly affected becomes progressively narrower. The comparison for different values of $\varepsilon$ therefore illustrates the response of the regularising net $(u_\varepsilon)_\varepsilon$ to the spatial concentration of the singular diffusion coefficient. Away from the singular point, the profiles are considerably less sensitive to further concentration of the mollifier.

\begin{figure}[h]
\includegraphics[scale=0.3]{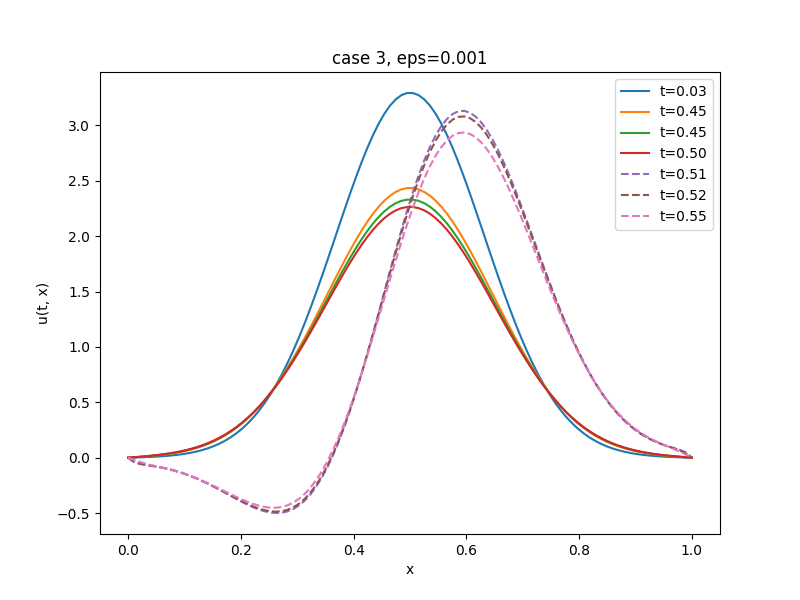}
\includegraphics[scale=0.3]{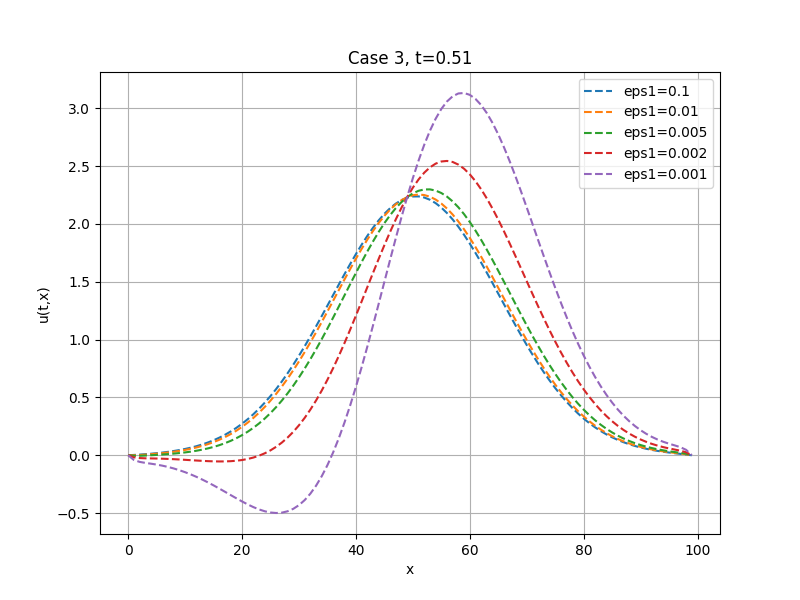}
\caption{Case 3}
\label{case3-1}
\end{figure}

Figure \ref{case3-1} illustrates the behaviour of the numerical solution $u(t,x)$, which serves as the solution to equation \ref{Equation_num2}, analysed for $\varepsilon$=0.001 in Case 3. In this scenario, the coefficient $b(t,x)$ is represented by a Dirac delta function supported at point $t=0.5$. The right panel compares the solution profiles for various values $\varepsilon=0.1$, $0.01$, $0.005$, $0.002$, and $0.0001$.

Figure \ref{case3-1} shows the evolution of the heat distribution, which decays smoothly until the singularity time $t = 0.25$. At this instant, the solution exhibits a discontinuous jump to the right, a phenomenon driven by the advective term $b(t, x)u_x(t, x)$. Following the jump, the distribution resumes its smooth decay.

As $\varepsilon$ decreases, the support of the regularised Dirac distribution in $b_\varepsilon$ becomes narrower, while its amplitude increases in such a way that its total mass remains fixed. Consequently, the influence of the singular drift is concentrated into an increasingly short time interval around $t=0.5$. The right panel illustrates this increasingly rapid transition of the solution. Thus, the relevant effect as $\varepsilon\to0$ is the temporal concentration of the action of the drift rather than a simple inverse proportionality between the magnitude of the solution displacement and $\varepsilon$.

\begin{figure}[h]
\includegraphics[scale=0.3]{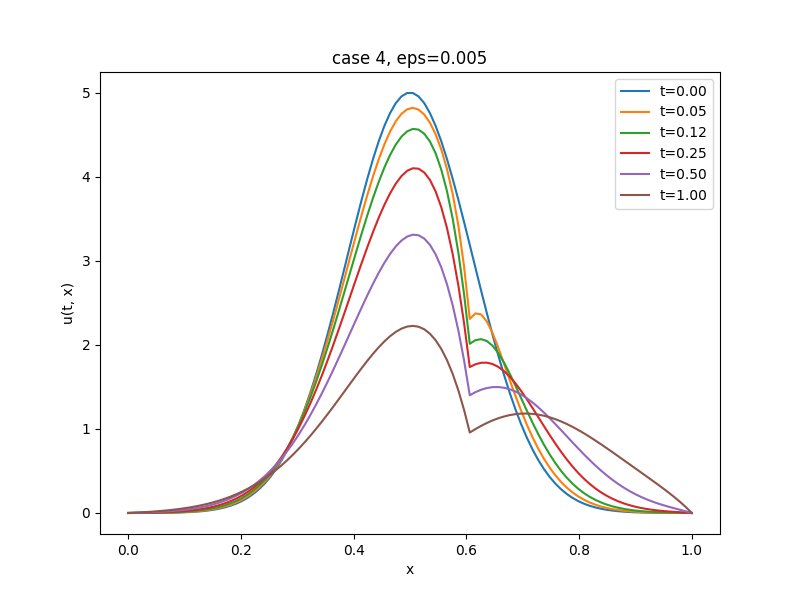}
\includegraphics[scale=0.3]{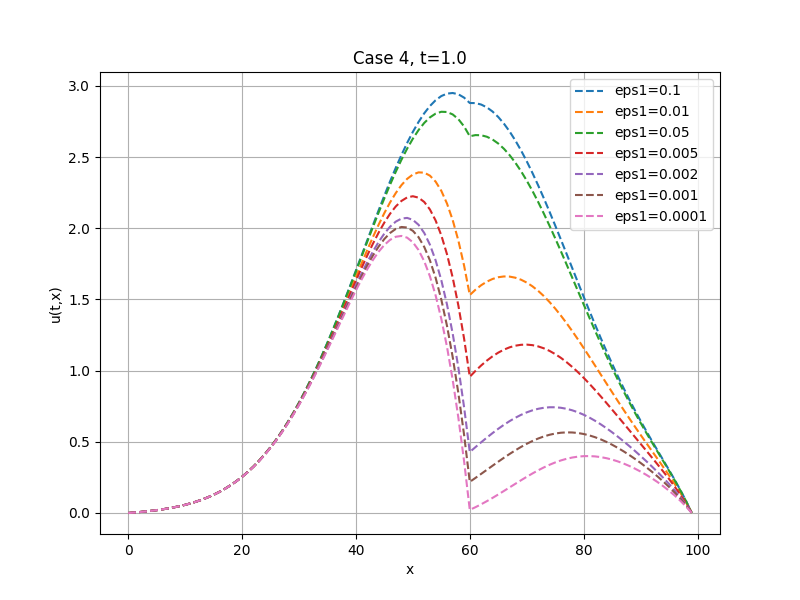}
\caption{Case 4}
\label{case4-1}
\end{figure}

Figure \ref{case4-1} analyses the heat distribution under a Dirac delta potential $q(x)$ at $x=0.6$. On the left, the evolution of the distribution is shown for a fixed $\varepsilon=0.005$. On the right, the results at the final time are compared for different values of $\varepsilon$.

Figure \ref{case4-1} shows that the overall heat distribution decays over time. Concurrently, the Dirac delta potential at $x=0.6$ induces a localised cooling effect. As illustrated on the left, this localised cooling effect occurs more rapidly for smaller values of $\varepsilon$.

\begin{figure}[h]
\includegraphics[scale=0.3]{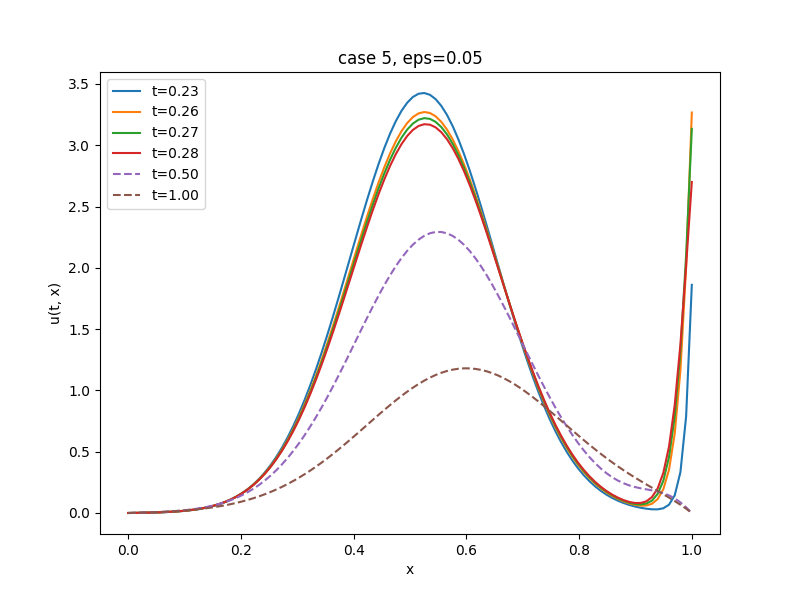}
\includegraphics[scale=0.3]{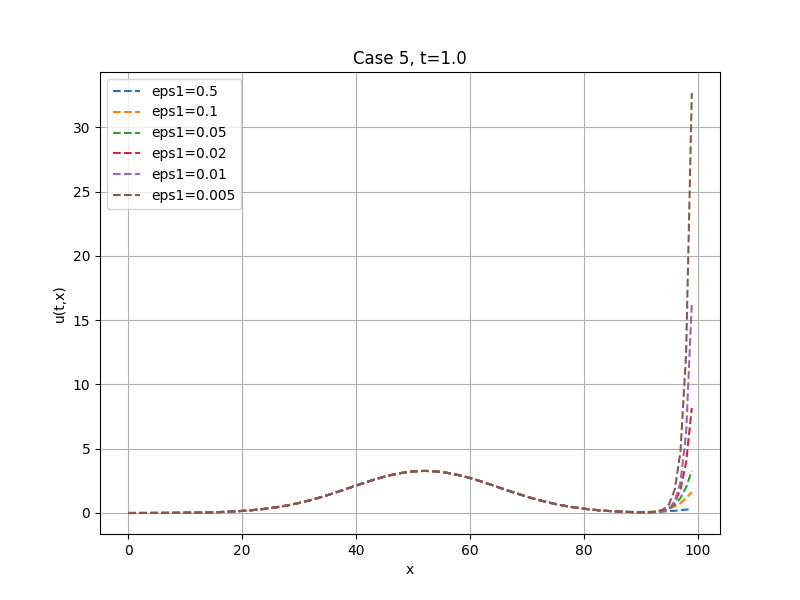}
\caption{Case 5}
\label{case5-1}
\end{figure}

Figure \ref{case5-1} illustrates the heat distribution subject to an inhomogeneous Dirichlet boundary condition with $g_1(t)$ modelled by a Dirac delta function at $t=0.25$. The left figure shows the distribution's evolution for a fixed $\varepsilon=0.005$. The right figure compares the profiles at the final time for various values of $\varepsilon$.

As in previous cases, the heat distribution in Figure \ref{case5-1} decays over time. At $t=0.25$, the Dirac delta boundary condition activates, causing a localised heating effect near the right boundary. This temperature increase is transient; after several time steps, the boundary temperature decreases. The right panel shows that this localised heating occurs more rapidly for smaller values of $\varepsilon$.

As $\varepsilon$ decreases, the boundary excitation becomes larger in amplitude but acts over an increasingly short time interval, while its total mass remains fixed. Hence, the influence of the singular boundary datum becomes increasingly concentrated near the singular time. The corresponding effect first appears close to the right boundary and is then smoothed and propagated into the interior of the spatial domain.

\subsection{Two-dimensional simulations}

We complement the preceding one-dimensional experiments with two-dimensional
simulations for Cases~2 and~4. Case~1 is regular and introduces no singular spatial effect, whereas the singularity in Case~3 is purely time-dependent; hence, their two-dimensional extensions would not provide additional qualitative information. Let
$
\Omega=(0,1)^2,
\,\,
\mathbf{x}=(x,y),
$
and consider the regularised initial--boundary value problem
\begin{align}\label{Equation_num_2D}
\partial_t u_\varepsilon(t,\mathbf{x})
&-
\nabla\cdot
\left(
a_\varepsilon(\mathbf{x})
\nabla u_\varepsilon(t,\mathbf{x})
\right)
+
\mathbf{b}\cdot\nabla u_\varepsilon(t,\mathbf{x})
\nonumber\\
&+
q_\varepsilon(\mathbf{x})
u_\varepsilon(t,\mathbf{x})
=
0,
\qquad
(t,\mathbf{x})\in(0,T]\times\Omega,
\end{align}
subject to the homogeneous Dirichlet boundary condition
\[
u_\varepsilon(t,\mathbf{x})=0,
\qquad
(t,\mathbf{x})\in(0,T]\times\partial\Omega,
\]
and the initial condition
\[
u_\varepsilon(0,x,y)
=
5\exp\left(
-\frac{(x-0.5)^2+(y-0.5)^2}{0.025}
\right).
\]
The vector \(\mathbf{b}\) represents the two-dimensional transport coefficient.

For a Dirac distribution concentrated at the interior point
$
\mathbf{x}_*=(0.45,0.45),
$
we use the tensor-product regularisation
\[
\delta_{\varepsilon,\mathbf{x}_*}(x,y)
=
\frac{1}{\varepsilon^2}
\varphi\left(\frac{x-0.45}{\varepsilon}\right)
\varphi\left(\frac{y-0.45}{\varepsilon}\right).
\]
Its total mass is equal to one, while its support becomes increasingly
concentrated around \(\mathbf{x}_*\) as \(\varepsilon\to0\).

The uniform grid discretises the spatial domain
$
x_i=ih_x,
\,\,
y_j=jh_y,
$
and the time derivative is approximated by the backward Euler method. The
two-dimensional diffusion operator is approximated using the standard
five-point finite-difference stencil, while the first-order transport terms
are approximated by one-sided differences. Thus, the two-dimensional
computations use the same semi-implicit principle as the one-dimensional
scheme described above. At every time step, the resulting sparse linear
system is solved for the grid values
$
u_{i,j}^{n}
\approx
u_\varepsilon(t_n,x_i,y_j).
$

In the two-dimensional version of Case~2, the singularity occurs in the
diffusion coefficient:
\[
a_\varepsilon(x,y)
=
1+\delta_{\varepsilon,\mathbf{x}_*}(x,y),
\qquad
q_\varepsilon(x,y)=1.
\]
In the two-dimensional version of Case~4, the diffusion coefficient is
regular, whereas the singularity occurs in the potential:
\[
a_\varepsilon(x,y)=1,
\qquad
q_\varepsilon(x,y)
=
1+\delta_{\varepsilon,\mathbf{x}_*}(x,y).
\]

\begin{figure}[h]
\centering
\includegraphics[width=\textwidth]{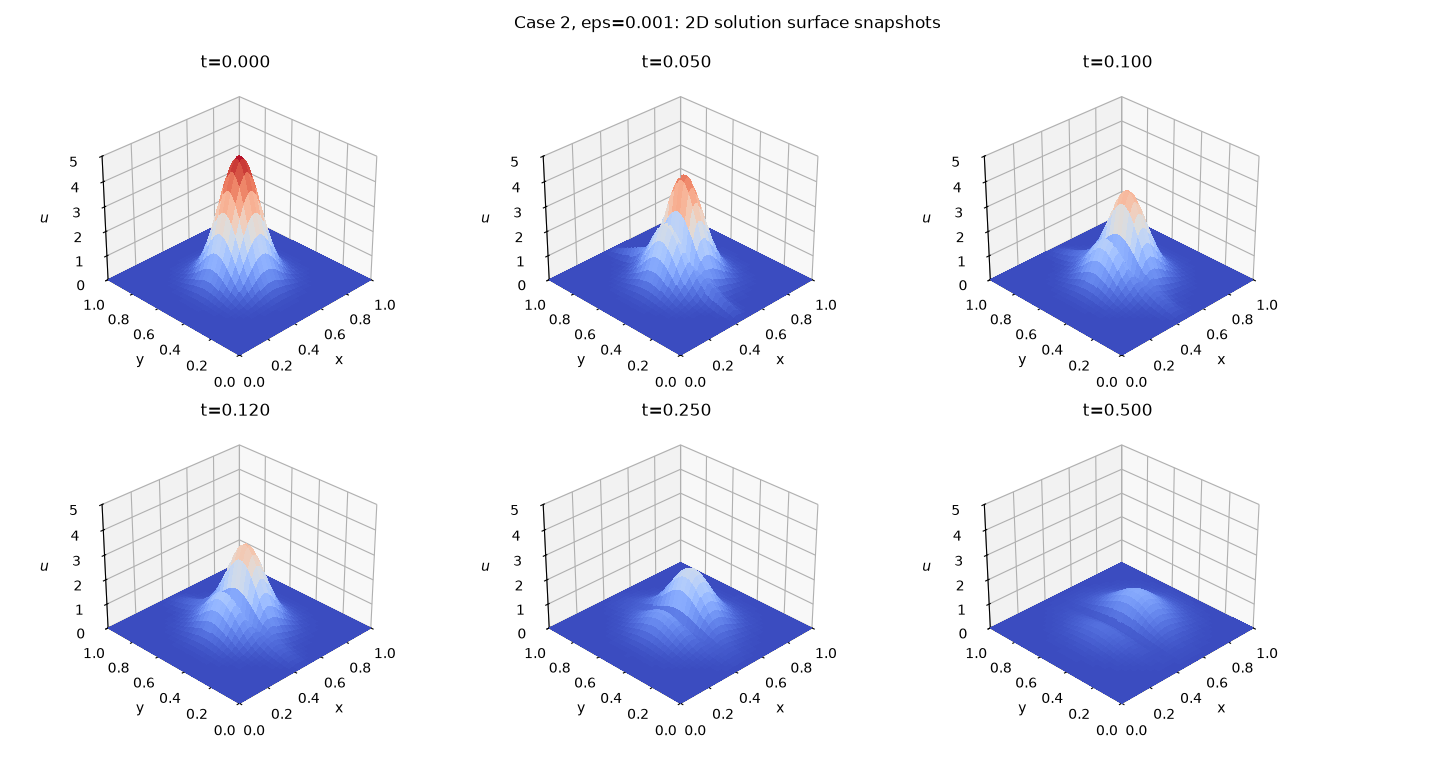}
\caption{Two-dimensional numerical solution for Case~2 with
\(\varepsilon=0.001\) at different time levels.}
\label{case2-2D}
\end{figure}

Figure~\ref{case2-2D} presents surface snapshots of the numerical solution for
the two-dimensional version of Case~2 with \(\varepsilon=0.001\).

The initial heat distribution is concentrated near the centre of the domain.
As time increases, the peak decreases because of diffusion and the positive
potential term, while the transport term produces a displacement of the
distribution. The singular contribution to \(a_\varepsilon\) locally
increases the diffusivity near \(\mathbf{x}_*\). Consequently, the solution
is smoothed more strongly in a neighbourhood of the singular point, and the
initial peak becomes progressively lower and broader. The effect is localised
near the support of the regularised Dirac distribution, while the solution
continues to evolve smoothly in the remaining part of the domain.

\begin{figure}[h]
\centering
\includegraphics[width=\textwidth]{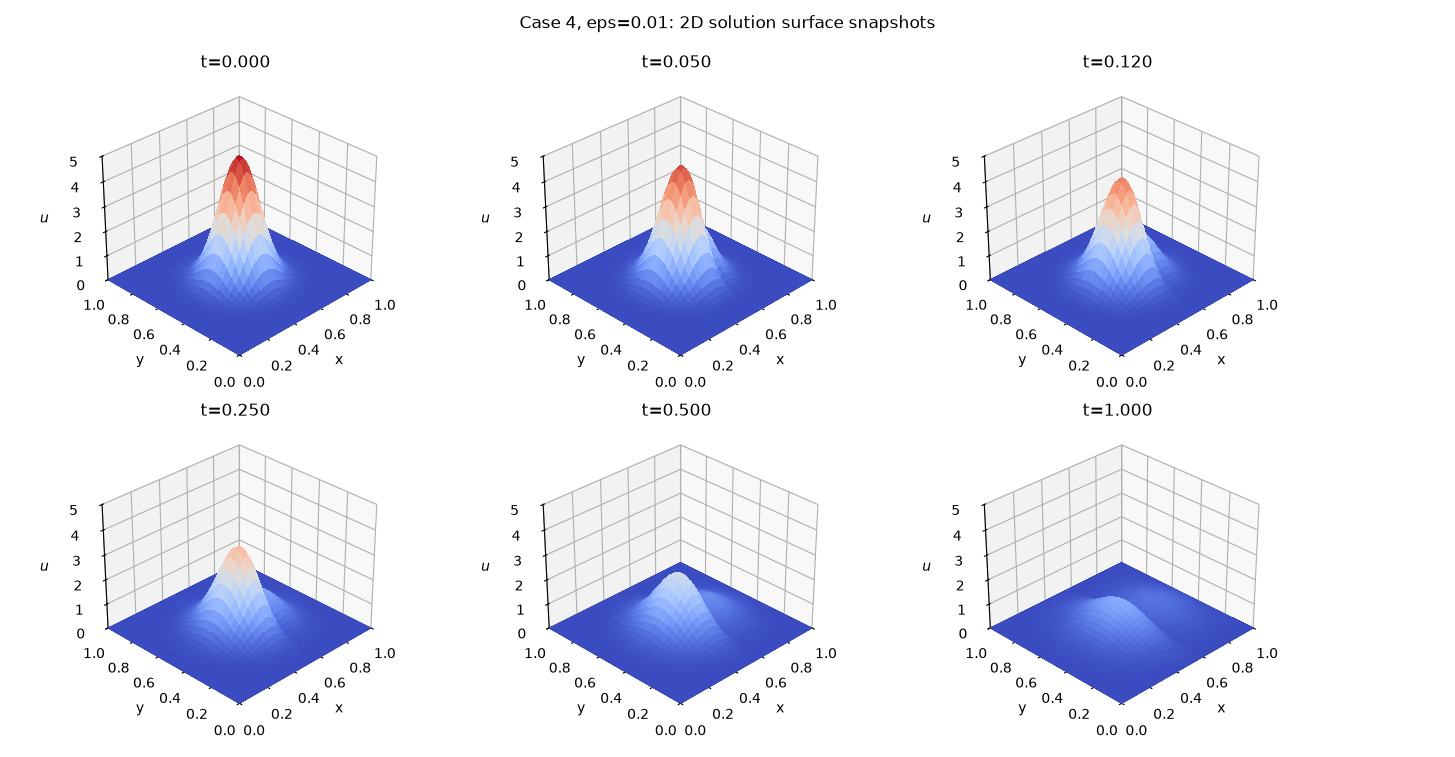}
\caption{Two-dimensional numerical solution for Case~4 with
\(\varepsilon=0.01\) at different time levels.}
\label{case4-2D}
\end{figure}

Figure~\ref{case4-2D} shows the corresponding numerical solution for the
two-dimensional version of Case~4 with \(\varepsilon=0.01\).

In this case, the regularised Dirac distribution occurs in the potential
coefficient. Since
$
q_\varepsilon(x,y)
=
1+\delta_{\varepsilon,\mathbf{x}_*}(x,y)
\geq 1,
$
it acts as a localised damping term. The surface plots show the decay of the
solution amplitude over time together with the displacement caused by the
transport term. Near \(\mathbf{x}_*\), the singular potential produces
additional local damping, whereas diffusion propagates and smooths this
effect throughout the domain.

The two-dimensional simulations confirm that the qualitative mechanisms
observed in the one-dimensional experiments persist in higher dimensions.
A singular diffusion coefficient produces localised enhanced smoothing,
whereas a positive singular potential produces localised damping. These
figures are intended to illustrate the spatial structure of the regularised
solutions for representative values of \(\varepsilon\). The dependence of
the regularising net \((u_\varepsilon)_\varepsilon\) on decreasing
\(\varepsilon\) is demonstrated primarily by the one-dimensional comparisons
presented above.

Overall, the numerical experiments illustrate how the regularising net $(u_\varepsilon)*\varepsilon$ responds when the mollified coefficients and data become increasingly concentrated as $\varepsilon\to0$. For spatially singular coefficients, the influence of the singularity becomes increasingly localised near its spatial support, whereas temporally singular coefficients and boundary data act over increasingly short time intervals. Away from the corresponding singular regions, the numerical solution profiles are considerably less sensitive to further decreases of $\varepsilon$. These observations are consistent with the role of regularisation in the construction of very weak solutions: for each fixed $\varepsilon>0$, one solves a classical uniformly parabolic problem, while the family $(u*\varepsilon)*\varepsilon$ describes the behaviour associated with the original distributional coefficients and data. The numerical experiments are intended to illustrate this regularisation procedure and the qualitative behaviour of its representatives.


\begin{thebibliography}{99}

\bibitem{ARST21c}
Altybay, A., Ruzhansky, M., Sebih, M.E., Tokmagambetov, N.,
The heat equation with strongly singular potentials,
Appl. Math. Comput. 399 (2021) 126006.

\bibitem{CDRT23}
Chatzakou, M., Dasgupta, A., Ruzhansky, M., Tushir, A.,
Discrete heat equation with irregular thermal conductivity and
tempered distributional data,
Proc. R. Soc. Edinb. Sect. A Math. (2023) 1--24.

\bibitem{CRT22a}
Chatzakou, M., Ruzhansky, M., Tokmagambetov, N.,
The heat equation with singular potentials. II: Hypoelliptic case,
Acta Appl. Math. 179 (2022) 2.

\bibitem{DiNardoFeoGuibe13}
Di Nardo, R., Feo, F., Guibé, O.,
Uniqueness of renormalized solutions to nonlinear parabolic problems
with lower-order terms,
Proc. R. Soc. Edinb. Sect. A Math. 143 (2013), no.~6, 1185--1208.
\url{https://doi.org/10.1017/S0308210511001831}

\bibitem{Evans}
Evans, L.C.,
Partial Differential Equations,
American Mathematical Society, Providence, RI, 2010.

\bibitem{Friedlander}
Friedlander, F.G., Joshi, M.,
Introduction to the Theory of Distributions,
Cambridge University Press, Cambridge, 1998.

\bibitem{Gar21}
Garetto, C.,
On the wave equation with multiplicities and space-dependent
irregular coefficients,
Trans. Amer. Math. Soc. 374 (2021) 3131--3176.

\bibitem{GR15}
Garetto, C., Ruzhansky, M.,
Hyperbolic second-order equations with non-regular time-dependent
coefficients,
Arch. Ration. Mech. Anal. 217 (2015) 113--154.

\bibitem{Issoglio19}
Issoglio, E.,
A non-linear parabolic PDE with a distributional coefficient and its
applications to stochastic analysis,
J. Differential Equations 267 (2019), no.~10, 5976--6003.
\url{https://doi.org/10.1016/j.jde.2019.06.014}

\bibitem{LionsMagenes}
Lions, J.-L., Magenes, E.,
Non-Homogeneous Boundary Value Problems and Applications,
Vol.~I,
Springer-Verlag, Berlin, 1972.

\bibitem{MichelVovelle03}
Michel, A., Vovelle, J.,
Entropy formulation for parabolic degenerate equations with general
Dirichlet boundary conditions and application to the convergence of
finite-volume methods,
SIAM J. Numer. Anal. 41 (2003), no.~6, 2262--2293.
\url{https://doi.org/10.1137/S0036142902406612}

\bibitem{MM2005}
Morton, K.W., Mayers, D.F.,
Numerical Solution of Partial Differential Equations: An Introduction,
2nd ed., Cambridge University Press, Cambridge, 2005.

\bibitem{Ober}
Oberguggenberger, M.,
Multiplication of Distributions and Applications to Partial Differential Equations,
Pitman Research Notes in Mathematics Series, Vol.~259,
Longman Scientific \& Technical, Harlow, 1992.

\bibitem{PierreVovelle12}
Pierre, M., Vovelle, J.,
A kinetic approach in nonlinear parabolic problems with
\(L^1\)-data,
Z. Anal. Anwend. 31 (2012), no.~3, 307--334.
\url{https://doi.org/10.4171/ZAA/1462}

\bibitem{RSY22}
Ruzhansky, M., Shaimardan, S., Yeskermessuly, A.,
Wave equation for Sturm--Liouville operator with singular potentials,
J. Math. Anal. Appl. 531 (2024) 127783.
\url{https://doi.org/10.1016/j.jmaa.2023.127783}

\bibitem{RT17a}
Ruzhansky, M., Tokmagambetov, N.,
Very weak solutions of the wave equation for the Landau Hamiltonian
with irregular electromagnetic field,
Lett. Math. Phys. 107 (2017) 591--618.

\bibitem{RT17b}
Ruzhansky, M., Tokmagambetov, N.,
Wave equation for operators with discrete spectrum and irregular
propagation speed,
Arch. Ration. Mech. Anal. 226 (2017) 1161--1207.


\bibitem{RY22}
Ruzhansky, M., Yeskermessuly, A.,
Wave equation for Sturm--Liouville operator with singular intermediate
coefficient and potential,
Bull. Malays. Math. Sci. Soc. 46 (2023) 195.
\url{https://doi.org/10.1007/s40840-023-01587-y}

\bibitem{ARST25}
A. Altybay, M. Ruzhansky, M.E. Sebih, N. Tokmagambetov, Singular hyperbolic type equations and tsunami propagation for irregular topographies, J. Math. Sci. 291 (2) (2025).

\bibitem{RY24a}
Ruzhansky, M., Yeskermessuly, A.,
Heat equation for Sturm--Liouville operator with singular propagation
and potential,
J. Appl. Anal. (2024).
\url{https://doi.org/10.1515/jaa-2023-0146}

\bibitem{RY25}
Ruzhansky, M., Yeskermessuly, A.,
Singular Klein--Gordon equation on a bounded domain,
J. Math. Sci. 291 (2025) 97--116.
\url{https://doi.org/10.1007/s10958-025-07782-5}

\bibitem{Sch54}
Schwartz, L.,
Sur l'impossibilit\'e de la multiplication des distributions,
C. R. Acad. Sci. Paris 239 (1954) 847--848.

\bibitem{Strikwerda2004}
Strikwerda, J.C.,
Finite Difference Schemes and Partial Differential Equations,
2nd ed., SIAM, Philadelphia, 2004.

\bibitem{Thomee06}
Thomée, V.,
Galerkin Finite Element Methods for Parabolic Problems,
2nd ed., Springer Series in Computational Mathematics, Vol.~25,
Springer, Berlin, 2006.

\end{thebibliography}
\end{document}